\documentclass[leqno,a4paper,11pt]{article}
\usepackage{draheniostyle}

\usepackage[
pdfauthor={Xabier Legaspi},
pdftitle={Growth of quasi-convex subgroups in groups with a constricting element},
pdfkeywords={Exponential growth, hyperbolic groups, contracting elements, convex-cocompactness},
pdftex, 
breaklinks
]{hyperref}

\hypersetup{
	colorlinks=true,%
	citecolor=VerdeVento,%
	filecolor=VerdeVento,%
	linkcolor=VerdeVento,%
	urlcolor=VerdeVento,%
}

\addto\extrasenglish{%
}

\title{\normalsize\bfseries\color{VerdeVento} GROWTH OF QUASI-CONVEX SUBGROUPS\\ IN GROUPS WITH A CONSTRICTING ELEMENT}
\author{\small\bfseries by Xabier LEGASPI}
\date{\small 20th June 2023}

\begin{document}
	
	\maketitle
	
	% Abstract
	\begin{center}
		\rule{2cm}{0.6pt}
		\vspace{0.5cm}
	\end{center}
	
	\begingroup
	\centering\small
	\begin{minipage}{\dimexpr\paperwidth-10cm}
		\textsc{Abstract}. --- \footnotesize
		Given a group $G$ acting by isometries on a metric space $X$, we consider a preferred collection of paths of the space $X$, a \emph{path system}, and study the spectrum of relative exponential growth rates and quotient exponential growth rates of the infinite index subgroups of $G$ that are quasi-convex with respect to this path system. If $G$ contains a constricting element with respect to the same path system, we are able to determine when the growth rates of the first kind are strictly smaller than the growth rate of $G$, and when the growth rates of the second kind coincide with the growth rate of $G$. Examples of applications include relatively hyperbolic groups, $\CAT(0)$ groups and hierarchically hyperbolic groups containing a Morse element.
	\end{minipage}
	\par
	\endgroup
	\bigskip\noindent
	
	% Keywords and Math Subject Classification
	{
		\footnotesize
		
		\medskip
		
	}
	
	% Index
	\tableofcontents
	
	% Contents
	\section{Introduction}
\label{sec:introduction}

The action of a group $G$ on a metric space $X$ is called \emph{proper} if for every $r\ge 0$, and for every $x\in X$, the number of elements $u \in G$ moving $x$ at distance at most $r$ is finite. Let $G$ be a group acting properly by isometries on a metric space $X$. The \emph{relative exponential growth rate} of the action of a subset $U\subset G$ on $X$ is the number
$$\omega(U,X)=\limsup_{r\to \infty}\frac{1}{r}\log|\setc{u \in U}{|u o-o|\le r}|,$$
whose value is independent of the point $o\in X$. Let $H$ be a subgroup of $G$. Let $H_L$ and $H_R$ be respectively minimal left and right transversals of $H$ at $o$, i.e., such that for every $u\in H_L$ and $v\in H_R$, 
$$|uo-o|=\inf_{h\in H} |uho-o|,\quad \text{and } |vo-o|=\inf_{h\in H} |hvo-o|.$$ 
In this article we study the numbers
$$\omega(H):=\omega(H,X), \quad \omega(G/H):=\omega(H_L,X), \quad \text{and }\omega(H\backslash G):=\omega(H_R,X).$$
The values of $\omega(G/H)$ and $\omega(H\backslash G)$ do not depend on the choice of the minimal transversal. Consider the following general problem. When do $G$ and $H$ determine a solution to to the system of equations below?
\begin{align*}
	\begin{cases}
		\omega(H)<\omega(G),\\
		\omega(G/H)=\omega(G),\\
		\omega(H\backslash G)=\omega(G).
	\end{cases}
\end{align*}
We see from the definitions that 
$$\omega(H/G)=\omega(H\backslash G), \quad \text{and } 0 \le \max\qty{\omega(H), \omega(G/H)}\le \omega(G).$$
In the extreme case in which $H$ has finite index in $G$, one can easily prove that
\begin{align*}
	\begin{cases}
		\omega(H)=\omega(G),\\
		\omega(G/H)=0.
	\end{cases}
\end{align*}
In general, it is a hard problem to obtain precise estimations of relative exponential growth rates of infinite index subgroups. However, it is known, \cite{dahmani_growth_2019, antolin_counting_2021, gitik_growth_2020}, that if $G$ is a non-virtually cyclic group acting geometrically on a hyperbolic space $X$ and $H$ is an infinite index quasi-convex subgroup of $G$, then
\begin{align*}
	\begin{cases}
		\omega(H)<\omega(G),\\
		\omega(G/H)=\omega(G).
	\end{cases}
\end{align*}
The arguments of \cite{dahmani_growth_2019, antolin_counting_2021} are based on automatic structures and regular languages, with influence of the works of J. Cannon \cite{cannon_combinatorial_1984,cannon_theory_1991}. This fact also influenced other authors that partially extended the hyperbolic case result, \cite{cordes_regularity_2022}. In Chapter $1$ we go beyond the hyperbolic case and we obtain two main results (\autoref{th:entropy1} and \autoref{th:entropy2}) with elementary proofs that do not require the theory of regular languages and automata. We will be interested in groups acting properly on metric spaces conditioned by a very general notion of ``non-positive curvature'' introduced by A. Sisto in \cite{sisto_contracting_2018} --- \emph{containing a constricting element with respect to a path system} --- while the infinite index subgroups object of our study will satisfy a very general notion of ``convex cocompactness'' --- \emph{quasi-convexity with respect to a path system}. 

The remaining of this section is structured as follows. First of all, we will mention two applications. Later we will give an informal explanation of our general setting as the result of a natural generalisation of these applications. We expect that this will be enough to understand our main theorems stated right after that. We will give another application at the end.

\paragraph{Groups acting properly with a strongly contracting element.} Members of this class contain elements that ``behave like'' a loxodromic isometry in a hyperbolic space -- in a strong sense. Let $\delta\ge 0$. A \emph{subset} $A$ of $X$ is $\delta$-\emph{strongly contracting} if the diameter of the nearest-point projection on $A$ of any metric ball of $X$ not intersecting $A$ is less than $\delta$. An \emph{element} $g$ of $G$ is $\delta$-\emph{strongly contracting} if it has infinite order and there exists an orbit of the cyclic subgroup generated by $g$ that is $\delta$-strongly contracting. In his seminal paper M. Gromov introduced the concept of $\delta$-hyperbolic space, \cite{gromov_hyperbolic_1987}. He observed that most of the large scale features of negative curvature can be described in terms of thin triangles. Nowadays, there are plenty of reformulations of the $\delta$-hyperbolicity. In particular, H. Masur and Y. Minsky gave one by describing geodesics in terms of strong contraction:

\begin{ex}
	A geodesic metric space $X$ is hyperbolic if and only if there exists $\delta\ge 0$ such that any geodesic segment of $X$ is $\delta$-strongly contracting, \cite[Theorem~2.3]{masur_geometry_1999}.
\end{ex}

The following are some subclasses of groups acting properly with a strongly contracting element:

\begin{enumerate}[label=(\roman*)]
	\item $\mathbf{H}=$ ``$G$ is a group acting properly  with a loxodromic element on a hyperbolic space $X$.'' In $\mathbf{H}$, an element is loxodromic if and only if it is strongly contracting. See \cite{coornaert_geometrie_1990}.
	\item  $\mathbf{RH}=$ ``$G$ is a relatively hyperbolic group acting with a hyperbolic element on a locally finite Cayley graph $X$ of $G$.'' In $\mathbf{RH}$, hyperbolic elements are strongly contracting. See \cite[Corollary~1.7]{osin_elementary_2006} and \cite[Theorem~2.14]{sisto_projections_2013}.
	\item  $\mathbf{CAT_0}=$ ``$G$ is a group acting properly with a rank-one element on a proper $\CAT(0)$ space $X$.'' In $\mathbf{CAT_0}$, rank-one elements are strongly contracting. See \cite[Theorem~5.4]{bestvina_characterization_2009} and \cite{cashen_morse_2020}.
	\item  $\mathbf{Mod_T}=$ ``$G$ is the mapping class group of an orientable surface of genus $g$ and $p$ marked points of complexity $3g+p-4>0$ acting on its Teichmüller space endowed with the Teichmüller metric.'' In $\mathbf{Mod_T}$, pseudo-Anosov elements are strongly contracting. See \cite{minsky_quasi-projections_1996} and \cite[Proposition~4.6]{masur_geometry_1999}.
	\item $\mathbf{GSC}=$ ``$G$ is an infinite graphical small cancellation group associated to a $Gr'(1/6)$-labeled graph with finite
	components labeled by a finite set $S$, acting on the Cayley graph $X$ of $G$ with respect to $S$.'' In $\mathbf{GSC}$, loxodromic WPD elements for the action of $G$ on the hyperbolic coned-off Cayley graph constructed by D. Gruber and A. Sisto in \cite{gruber_infinitely_2018} are strongly contracting. See \cite[Theorem~5.1]{arzhantseva_negative_2019}.
	\item $\mathbf{Gar}=$ ``$G$ is the quotient of a $\Delta$-pure Garside group of finite type by its center, acting  with a Morse element on the Cayley graph $X$ of $G$ with respect to the Garside generating set.'' In $\mathbf{Gar}$, Morse elements are strongly contracting. See \cite[Theorem~5.5]{calvez_morse_2021}.
	\item $\mathbf{Inj}=$ ``$G$ is a group acting properly with a Morse element on an injective metric space $X$.'' In $\mathbf{Inj}$, an element is Morse if and only if it is strongly contracting. See \cite{sisto_morse_2023}.
\end{enumerate}

An appropriate notion of convex cocompactness in this setting is just the usual quasi-convexity. Let $\eta\ge 0$. A \emph{subset $Y$} of $X$ is $\eta$-quasi-convex if any geodesic of $X$ with endpoints in $Y$ is contained in the $\eta$-neighbourhood of $Y$. A \emph{subgroup} $H$ of $G$ is $\eta$-\emph{quasi-convex} if there exists an orbit of $H$ that is $\eta$-quasi-convex.

Our theorem below generalises \cite[Theorem~4.8]{yang_statistically_2019} and \cite[Theorems~1.1~and~1.3]{dahmani_growth_2019}:

\begin{thm}
	\label{th:strong contraction}
	If $G$ is a non-virtually cyclic group acting properly with a strongly contracting element on a geodesic metric space $X$, and $H$ is an infinite index quasi-convex subgroup of $G$, then 
	\begin{align*}
		\begin{cases}
			\omega(H)<\omega(G),\\
			\omega(G/H)=\omega(G).
		\end{cases}
	\end{align*}
\end{thm}

\paragraph{Hierarchically hyperbolic groups.} Let $\operatorname{Mod}(\Sigma_{g,p})$ be the mapping class group of an orientable surface $\Sigma_{g,p}$ of genus $g$ and $p$ marked points of complexity $3g+p-4>0$. We would like to apply \autoref{th:strong contraction} to $\operatorname{Mod}(\Sigma_{g,p})$ with respect to the word metric. However, we do not know whether $\operatorname{Mod}(\Sigma_{g,p})$ acts with a strongly contracting element on any of its locally finite Cayley graphs or not. Maybe the candidates that come to mind are the pseudo-Anosov elements, and evidence suggests that not all of them are strongly contracting: K. Rafi and Y. Verberne constructed a generating set $U$ of $\operatorname{Mod}(\Sigma_{0,5})$ and a pseudo-Anosov element which is not strongly contracting for the action of $\operatorname{Mod}(\Sigma_{0,5})$ on the Cayley graph of $\operatorname{Mod}(\Sigma_{0,5})$ with respect to $U$, \cite[Theorem~1.3]{rafi_geodesics_2021}. We were able to avoid this setback by looking into the class of hierarchically hyperbolic groups, introduced by J. Behrstock, M. Hagen and A.Sisto in \cite{behrstock_hierarchically_2017, behrstock_hierarchically_2019} as a generalisation of the Masur and Minsky hierarchy machinery of mapping class groups. Below we provide some examples of hierarchically hyperbolic groups. The reader should note that the metric space where they act with a hierarchically hyperbolic structure is any of their locally finite Cayley graphs:
\begin{enumerate}[label=(\roman*)]
	\item Mapping class groups of finite type surfaces, \cite{behrstock_hierarchically_2019}.
	\item Right-angled Artin groups, \cite{behrstock_hierarchically_2017}.
	\item Right-angled Coxeter groups, \cite{behrstock_hierarchically_2017}.
	\item Fundamental groups of $3$-manifolds without NIL or SOL components, \cite{behrstock_hierarchically_2019}.
\end{enumerate}

Now consider the following notion of convex cocompactness. A \emph{subset} $Y$ of $X$ is \emph{Morse} if for every $\kappa\ge 1$, $\lambda\ge 0$, there exists $\sigma\ge 0$ such that any $(\kappa,l)$-quasi-geodesic of $X$ with endpoints in $Y$ is contained in the $\sigma$-neighbourhood of $Y$. A \emph{subgroup} $H$ of $G$ is \emph{Morse} if there exists an orbit of $H$ that is Morse. An \emph{element} $g$ of $G$ is \emph{Morse} if it has infinite order and the cyclic subgroup generated by $g$ is Morse.

We have obtained the next result, partially generalising \cite[Theorem~A]{cordes_regularity_2022}:

\begin{thm}
	\label{th:hierarchically}
	If $G$ is a non-virtually cyclic hierarchically hyperbolic group acting on a locally finite Cayley graph $X$ of $G$ with a Morse element, and $H$ is an infinite index Morse subgroup of $G$, then
	\begin{align*}
		\begin{cases}
			\omega(H)<\omega(G),\\
			\omega(G/H)=\omega(G).
		\end{cases}
	\end{align*}
\end{thm}

We know that pseudo-Anosov elements of mapping class groups are Morse with respect to any word metric, \cite{behrstock_asymptotic_2006}, and that the infinite index Morse subgroups of the mapping class group are precisely the convex cocompact subgroups in the sense of mapping class groups, \cite[Theorem~A]{kim_stable_2019}, which allows us to obtain a more concrete statement:

\begin{cor}
	\label{cor:mapping-class-groups}
	If $G$ is the mapping class group of a surface of genus $g$ and $p$ marked points such that $3g+p-4>0$ acting on a locally finite Cayley graph $X$ of $G$, and $H$ is a convex cocompact subgroup of $G$, then 
	\begin{align*}
		\begin{cases}
			\omega(H)<\omega(G),\\
			\omega(G/H)=\omega(G).
		\end{cases}
	\end{align*}
\end{cor}

\begin{rem}
	Under the hypothesis of the previous corollary, we remark that the inequality $\omega(H)<\omega(G)$ was also obtained independently in \cite[Corollary~C]{cordes_regularity_2022}.
\end{rem}

\paragraph{Main results.} Now that we gave the big picture, we will give a technical definition that encapsulates the classes discussed so far. In order to do so, we make two observations. On the one hand, the strong contraction property can be reformulated in the following way. A subset $A$ of $X$ is \emph{strongly contracting} if and only if any geodesic segment of $X$ joining any pair of points $x,y\in X$ whose projections $p$ and $q$ via a nearest-point projection are far away passes next to $p$ and $q$, \cite[Proposition~2.9]{arzhantseva_growth_2015}. On the other hand, mapping class groups -- or more generally, hierarchically hyperbolic groups -- come with hierarchy paths, a family of special quasi-geodesics encoding substantial information about the geometry of the space and easier to work with than the set of all (quasi-)geodesics. For these reasons, in order to define very general notions of non-positive curvature and convex cocompactness, we will be considering path systems, introduced by A. Sisto in \cite{sisto_contracting_2018}:
\begin{df}[Path system group]
	Let $\mu\ge 1$, $\nu\ge 0$. A $(\mu,\nu)$-\emph{path system group} $(G,X,\scrP)$ is a group $G$ acting properly on a geodesic metric space $X$ together with a $G$-invariant collection $\scrP$ of paths of $X$ satisfying:
	\begin{enumerate}[label=(PS\arabic*)]
		\item $\scrP$ is closed under taking subpaths.
		\item For every $x,y\in X$, there exists $\gamma\in\scrP$ joining $x$ to $y$.
		\item Every element of $\scrP$ is a $(\mu,\nu)$-quasi-geodesic.
	\end{enumerate}
	We refer to $\scrP$ as $(\mu,\nu)$-\emph{path system}.
\end{df}
We fix $\mu\ge 1$, $\nu\ge 0$ and a $(\mu,\nu)$-path system group $(G,X, \scrP)$ for the following definitions. Let $\delta\ge 0$. We say that a subset $A$ of $X$ is $\delta$-\emph{constricting} if there exist a \emph{coarse nearest-point projection} of $X$ on $A$ with the property that any $\gamma\in \scrP$ joining any two pair of points $x,y\in X$ whose projections $p$ and $q$ are $\delta$-far away passes through the $\delta$-neighbourhoods of $p$ and $q$ (\autoref{def:constricting element path system}). An \emph{element} $g$ of $G$ is $\delta$-\emph{constricting} if it has infinite order and there exists a $\delta$-constricting orbit of the cyclic subgroup generated by $g$. Let $\eta\ge 0$. A subgroup $Y$ of $X$ is $\eta$-\emph{quasi-convex} if any $\gamma\in \scrP$ with endpoints in $Y$ is contained in the $\eta$-neighbourhood of $Y$ (\autoref{def:quasi-convexity path system}). A subgroup $H$ of $G$ is $\eta$-quasi-convex if there exist an $\eta$-quasi-convex orbit of $H$.

\begin{ex}
	\label{ex:applications}
	
	\begin{enumerate}[label=(\roman*)]
		\item Assume that the metric space $X$ is geodesic. An infinite order element of $G$ is strongly contracting if and only if it is constricting with respect to the set of all the geodesic segments of $X$, \cite[Proposition~2.9]{arzhantseva_growth_2015}.
		\item Assume that the group $G$ is hierarchically hyperbolic. An infinite order element $g$ of $G$ is Morse if and only if for every $\kappa\ge 1$, there exists $\delta\ge 0$ such that $g$ is $\delta$-constricting with respect to the set of all the $\kappa$-hierarchy paths. See \cite[Theorem~E]{russell_convexity_2021} and \cite[Lemma~1.27]{behrstock_quasiflats_2020}. 
	\end{enumerate}
\end{ex}

Finally, we state the main results of Chapter $1$. \autoref{th:strong contraction} and \autoref{th:hierarchically} are special cases. Our first result generalises work of W. Yang, \cite[Theorem~4.8]{yang_statistically_2019}, and F. Dahmani - D. Futer - D. Wise, \cite[Theorems~1.1~and~1.3]{dahmani_growth_2019}. The \emph{Poincaré series} $\mathscr{P}_{U}(s)$ based at $o\in X$ of a subset $U$ of $G$ is defined as
$$\forall\,s\ge 0,\quad \mathscr{P}_{U}(s)=\sum_{u\in U} e^{-s|uo-o|}$$
and modifies its behaviour at the relative exponential growth rate $\omega(U,X)$: the series diverges if $s<\omega(U,X)$ and converges if $s>\omega(U,X)$. At $s=\omega(U,X)$ the series can converge or diverge depending on the nature of $U$. This behaviour is independent of the point $o\in X$. We say that the action of $U$ on $X$ is \emph{divergent} if $\mathscr{P}_{U}(s)$ diverges at $s=\omega(U,X)$.

\begin{thm}[\autoref{th:entropy1-body}]
	\label{th:entropy1}
	Let $(G,X,\scrP)$ be a path system group. Assume that $G$ contains a constricting element. Let $H$ be an infinite index subgroup of $G$ satisfying the following:
	\begin{enumerate}[label=(\roman*)]
		\item $\omega(H)<\infty$.
		\item The action of $H$ on $X$ is divergent.
		\item $H$ is quasi-convex.
	\end{enumerate}	
	Then $\omega(H)<\omega(G)$.
\end{thm}

\begin{rem}
	Under the hypothesis of \autoref{th:entropy1}, one may ask if there is a growth gap, i.e, if
$$\sup_H \omega(H)<\omega(G),$$
where the supremum is taken among the infinite index subgroups $H$ of $G$ satisfying (i), (ii) and (iii). In our context, the answer is yes: there is a growth gap when $G$ is a hyperbolic group with \emph{Kazhdan's Property} (T), \cite[Theorem~1.2]{coulon_growth_2018}. However, one can show that there is no growth gap among free groups, \cite[Theorem~9.4]{dahmani_growth_2019}, or fundamental groups of compact special cube complexes,  \cite[Theorem~1.5]{li_no_2020}. The answer to our context could be different if one studied semigroups instead of subgroups, \cite[Theorem~A]{yang_statistically_2019}.
\end{rem}

In \cite[5.3.C]{gromov_hyperbolic_1987}, M. Gromov stated that in a torsion-free hyperbolic group $G$, any infinite index quasi-convex subgroup $H$ is a free factor of a larger quasi-convex subgroup. Gromov's ideas were later developed by G. N. Arzhantseva in \cite[Theorem~1]{arzhantseva_quasiconvex_2001}. More recently, J. Russell, D. Spriano and H. C. Tran generalised her result to the context of groups with the ``Morse local-to-global property'',  \cite[Corollary~3.5]{russell_local--global_2022}. Further, the problem seems connected to the ``$P_{\emph{Naive}}$ property'' studied by C. Abbott and F. Dahmani in the context of groups acting acylindrically on a hyperbolic space, \cite{abbott_property_2019}. In our context, we have obtained the following, in which there is no torsion-free assumption. We will see that \autoref{th:entropy1} is, in part, a consequence of this result:

\begin{thm}[\autoref{prop:injectivity}]
	\label{th:amalgam}
	Let $(G,X,\scrP)$ be a path system group. Assume that $G$ contains a constricting element $g_0$. Let $H$ be an infinite index quasi-convex subgroup of $G$. Then, there exist an element $g\in G$ conjugate to a large power of $g_0$ and a finite extension $E$ of $\langle g\rangle$ such that the intersection $H\cap E$ is finite and the natural morphism $H\ast_{H\cap E} \langle g, H\cap E\rangle\to G$ is injective.
\end{thm}

According to \autoref{prop:elementary properties} (6), the subgroup generated by a constricting element is always Morse, and in particular quasi-convex. Hence \autoref{th:amalgam}, for the choice of $H= \group{g_0}$, implies the following weak Tits alternative:

\begin{cor}
	\label{cor:Tits alternative}
	Let $(G,X,\scrP)$ be a path system group. Assume that $G$ contains a constricting element. Then, either $G$ is virtually cyclic or contains a free subgroup of rank two.
\end{cor}
\begin{rem}
To the best of our knowledge, the previous corollary has not been recorded for the class of groups acting properly with a strongly contracting element. The Tits alternative is known for hierarchically hyperbolic groups \cite[Theorem~9.15]{durham_boundaries_2017}, which is a much stronger result.
\end{rem}

In our second result we generalise work of Y. Antolín, \cite[Theorem~3]{antolin_counting_2021}, and R. Gitik - E. Rips, \cite[Theorem~2]{gitik_growth_2020}:

\begin{thm}
	\label{th:entropy2}
	Let $(G,X,\scrP)$ be a path system group. Assume that $G$ contains a constricting element. Let $H$ be an infinite index quasi-convex subgroup of $G$. Then 
	$$\omega(G/H)=\omega(G).$$
\end{thm}
Note that the study of \cite[Theorem~2]{gitik_growth_2020} concerns double cosets in the hyperbolic group case. We remark that in \cite[VII~D~39]{de_la_harpe_topics_2000}, P. de la Harpe says about the growth of double cosets: ``this theme has not received yet too much attention, but probably should''. In our context, for sake of simplicity, we decided to study single cosets instead, but one could possibly extend our result. Further, we remark that our result is connected to the study of I. Kapovich on the hyperbolicity and amenability of the Schreier graphs of infinite index quasi-convex subgroups of hyperbolic groups, \cite{kapovich_geometry_2002, kapovich_nonamenability_2002}. There's also work of A. Vonseel concerning the number of ends, \cite{vonseel_ends_2018}.

\begin{rem}
	\label{rem:infty}
	\begin{enumerate}[label=(\roman*)]
	\item Our main results \autoref{th:entropy1} and \autoref{th:entropy2} hold in the case $\omega(G)=\infty$. For instance, if $G$ is a group acting properly on a metric space $(X,|\cdot|)$, then we can define a new metric $|\cdot|'$ on $X$ by
	$$\forall\,x,y\in X,\quad |x-y|'=e^{-|x-y|}\cdot |x-y|.$$
	The metric distorts the growth of the orbit of $G$ exponentially. If $\omega(G)>0$ with respect to $|\cdot|$, then $\omega(G)=\infty$ with respect to $|\cdot|'$.
	\item If $G$ is a group acting geometrically on a metric space $X$, then $\omega(G)<\infty$.
	\end{enumerate}
\end{rem}

Now we are going to record a joint corollary to \autoref{th:entropy1} and \autoref{th:entropy2}. In general, it is not easy to decide whether the action of a groups is divergent or not. However, the following is a well-known consequence of \emph{Fekete's Subadditive Lemma}:

\begin{lem}[{\cite[Proposition~4.1~(1)]{dalbo_growth_2011}}]
	\label{lem:divergence}
	Let $G$ be a group acting properly on a geodesic metric space $X$. Let $o\in X$. Let $H\le G$ be a quasi-convex subgroup (in the classical sense). Then 
	$$\omega(H)=\inf_{n\ge 1} \frac{1}{n} \log| \setc{h\in H}{|ho-o|\le n}|=\lim_{n\to \infty}\frac{1}{n}\log |\setc{h\in H}{|ho-o|\le n}|.$$
	In particular $\omega(H)<\infty$. If in addition $H$ is infinite, then the action of $H$ on $X$ is divergent. 
\end{lem}

Combining \autoref{lem:divergence} with \autoref{cor:Tits alternative}, we obtain:

\begin{cor}
	\label{cor:applications}
	Let $(G,X,\scrP)$ be a path system group. Assume that $G$ is non-virtually cyclic and contains a constricting element.
	\begin{enumerate}[label=(\roman*)]
		\item If $\scrP$ is the set of all the geodesic segments of $X$, then for every infinite index quasi-convex subgroup $H$ of $G$, we have 
		\begin{equation*}
			\begin{cases}
				\omega(H)<\omega(G),\\
				\omega(G/H)=\omega(G).
			\end{cases}
		\end{equation*}
		\item For every infinite index Morse subgroup $H$ of $G$, we have
		\begin{equation*}
			\begin{cases}
				\omega(H)<\omega(G),\\
				\omega(G/H)=\omega(G).
			\end{cases}
		\end{equation*}
	\end{enumerate}
\end{cor}

\begin{rem}
	One can prove that the class of groups acting properly with a constricting element with respect to a path system is invariant under equivariant quasi-isometries. However, strongly contracting elements are not preserved under equivariant quasi-isometries,  \cite[Theorem~4.19]{arzhantseva_negative_2019}. In particular, \autoref{cor:applications} applies for instance to the action on a locally finite Cayley graph of any group acting geometrically on a $\CAT(0)$ space with a rank-one element.
\end{rem}

\begin{rem}
	The proofs of \autoref{th:strong contraction}, \autoref{th:hierarchically} and \autoref{cor:mapping-class-groups} now follow from our main results (\autoref{th:entropy1} and \autoref{th:entropy2}) in view of \autoref{ex:applications} and \autoref{rem:infty} (ii).
\end{rem}

\paragraph{Hierarchical quasi-convexity.} In hierarchically hyperbolic groups there is a notion of convex cocompactness more natural than Morseness. Let $G$ be a hierarchically hyperbolic group. A subgroup $H$ of $G$ is hierarchically quasi-convex if and only if for every $\kappa\ge 1$, there exists $\eta\ge 0$ such that $H$ is $\eta$-quasi-convex with respect to the set of all the $\kappa$-hierarchy paths of $G$, \cite[Proposition~5.7]{russell_convexity_2021}. Finally, in view of \autoref{rem:infty} (ii) and \autoref{ex:applications} (ii), we deduce two more applications from \autoref{th:entropy1} and \autoref{th:entropy2}:

\begin{thm}
	If $G$ is a hierarchically hyperbolic group acting on a locally finite Cayley graph $X$ of $G$ with a Morse element, and $H$ is an infinite index subgroup of $G$ satisfying:
	\begin{enumerate}[label=(\roman*)]
		\item the action of $H$ on $X$ is divergent,
		\item $H$ is hierarchically quasi-convex,
	\end{enumerate}	
	then $\omega(H)<\omega(G)$.
\end{thm}

\begin{thm}
	If $G$ is a hierarchically hyperbolic group acting on a locally finite Cayley graph $X$ of $G$ with a Morse element, and $H$ is an infinite index hierarchically quasi-convex subgroup of $G$, then $\omega(G/H)=\omega(G)$.
\end{thm}

\paragraph{Outline of the paper.}
In \autoref{sec:path-system} we will introduce the definitions of path system group, quasi-convex subgroup and constricting element. In \autoref{sec:growth-estimation} we will explain the two criteria that we will use to estimate the growth of quasi-convex subgroups. The rest of the chapter is devoted to the development of our geometric framework so that we can apply these criteria. In \autoref{sec:qc-int-img} we will prove a version of the \emph{bounded geodesic image property} of hyperbolic spaces. In \autoref{sec:buffering} we will introduce the notion of \emph{buffering sequence} and we will give a version of \emph{Behrstock inequality}. In \autoref{sec:finding}, given an infinite index quasi-convex subgroup and a quasi-convex element, we will produce another quasi-convex element whose orbit is ``transversal'' to the given subgroup. The proofs of both of our main results (\autoref{th:entropy1} and \autoref{th:entropy2}) share this argument. In \autoref{sec:constricting-elements} we will study the elementary closures of constricting elements apart from some geometric separation properties. Finally, in \autoref{sec:growth} we will prove our main results (including \autoref{th:amalgam}) by constructing an appropriate buffering sequence for each problem.

\paragraph{Acknowledgements.}
The author is deeply indebted to his supervisors Yago Antolín and Rémi Coulon. He would like to thank them for their many comments, suggestions, words of support and for reading earlier versions of this draft that helped to improve it from scratch. He would also like to acknowledge Adrien Le Boudec, Jacob Russell, Alessandro Sisto, Juan Souto, Markus Steenbock and the anonymous referee for their comments and discussions. The author is grateful to the Centre Henri Lebesgue ANR-11-LABX-0020-01 for creating an attractive mathematical environment and to the Institut Henri Poincaré (UAR 839 CNRS-Sorbonne Université) for its hospitality and support (through LabEx CARMIN, ANR-10-LABX-59-01) during the trimester program \emph{Groups Acting on Fractals, Hyperbolicity and Self-Similarity} in Spring 2022. The author was supported by grant Outgoing Mobility 2021 funded by Rennes Métropole; and grants SEV-2015-0554-18-4 and PID2021-126254NB-I00 funded by MCIN/AEI/10.13039/501100011033.

	\section{Path system geometry}
\label{sec:path-system}

This section is devoted to present the notations and vocabulary of the main geometric objects of this chapter. We formalise our notions of ``convex cocompactness'' and ``non-positive curvature''.

\paragraph{Metric geometry.} 

Let $X$ be a metric space. Given two points $x, x'\in X$, we write $|x-x'|$ for the distance between them. The \emph{ball of} $X$ of center $x\in X$ and radius $r\ge 0$ is
$$B_X(x,r)=\setc{y\in X}{|x-y|\le r}.$$
The \emph{distance between a point} $x\in X$ \emph{and a subset}  $Y\subset X$ is
$$d(x,Y)=\inf\,\setc{|x-y|}{y\in Y}.$$
Let $\eta\ge 0$. The $\eta$-\emph{neighbourhood} of a subset $Y\subset X$ is
$$Y^{+\eta}=\setc{x\in X}{d(x,Y)\le\eta}.$$
The \emph{distance} between two subsets $Y,Z\subset X$ is
$$d(Y,Z)=\inf\, \setc{|y-z|}{y\in Y, z\in Z}.$$  
The \emph{Hausdorff distance} between two subsets $Y,Z\subset X$ is
$$d_{\operatorname{Haus}}(Y,Z)=\inf\, \setc{\varepsilon\ge 0}{Y\subset Z^{+\varepsilon} \text{ and } Z\subset Y^{+\varepsilon}}.$$

\paragraph{Path system spaces.} 

Let $X$ be a metric space. A \textit{path} is a continuous map $\alpha \colon [a,b]\to X$. The \textit{initial and terminal points} of $\alpha$ are $\alpha(a)$ and $\alpha(b)$, respectively. They form the \textit{endpoints} of $\alpha$. We will frequently identify a path and its image. A \textit{subpath} of $\alpha$ is a restriction of $\alpha$ to a subinterval of $[a,b]$. The path $\alpha$ \textit{joins} the point $x\in X$ to the point $y\in X$ if $\alpha(a) = x$ and $\alpha(b) = y$. Note that for every $x,y\in\alpha$ there may be more than one subpath of $\alpha$ joining $x$ to $y$, unless the points are given by the parametrisation of $\alpha$. The \textit{length} of a path $\alpha$ is denoted by $\ell(\alpha)$. Unless otherwise stated a path is a \textit{rectifiable path} parametrised by \textit{arc length}. Let $\kappa\ge 1$, $l\ge 0$. A path $\alpha \colon [a,b]\to X$ is a $(\kappa,l)$-\textit{quasi-geodesic} if for every $t,t'\in [a,b]$,
$$|\alpha(t)-\alpha(t')|\le |t-t'|\le \kappa |\alpha(t)-\alpha(t')|+l.$$
Note that that $\ell(\alpha_{|[t,t']})=|t-t'|$. The following captures the idea of endowing a metric space with a collection of preferred paths.

\begin{df}[Path system space]	
	Let $\mu\ge 1$, $\nu\ge 0$. A $(\mu,\nu)$-\emph{path system space} $(X,\scrP)$ is a metric space $X$ together with a collection $\scrP$ of paths of $X$ satisfying:
	
	\begin{enumerate}[label=(PS\arabic*)]
		\item $\scrP$ is closed under taking subpaths.
		\item For every $x,y\in X$, there exists $\gamma\in\scrP$ joining $x$ to $y$.
		\item Every element of $\scrP$ is a $(\mu,\nu)$-quasi-geodesic.
	\end{enumerate}
	We refer to $\scrP$ as $(\mu,\nu)$-\emph{path system}.
\end{df}

We fix $\mu\ge 1$, $\nu\ge 0$ and a $(\mu,\nu)$-\emph{path system space} $(X,\scrP)$.

\begin{df}[Quasi-convex subset]
	Let $\eta\ge 0$. A \emph{subset} $Y\subset X$ is $\eta$-\emph{quasi-convex} if every $\gamma\in \scrP$ with endpoints in $Y$ is contained in the $\eta$-neighbourhood of $Y$.
\end{df}

\begin{df}[Constricting subset]
	\label{def:constricting}
	
	Let $\delta \ge 0$. A \emph{subset} $A\subset X$ is $\delta$-\emph{constricting} if there exists a map $\pi_A\colon X\to A$ satisfying:
	\begin{enumerate}[label=(CS\arabic*)]
		\item \textbf{Coarse retraction.}
		
		\noindent
		For every $x\in A$, we have $|x-\pi_A(x)|\le \delta$.
		
		\item  \textbf{Constriction.}
		
		\noindent
		For every $x,y\in X$ and for every $\gamma\in \scrP$ joining $x$ to $y$, if we have $|\pi_A(x)-\pi_A(y)|> \delta$, then $\gamma\cap B_X(\pi_A(x),\delta)\neq \varnothing$ and $\gamma\cap B_X(\pi_A(y),\delta)\neq \varnothing$.
	\end{enumerate}
	We refer to $\pi_A \colon X\to A$ as $\delta$-\emph{constricting map}.
\end{df}

\begin{figure}[htbp]
	\centering
	\def\svgwidth{\columnwidth}
	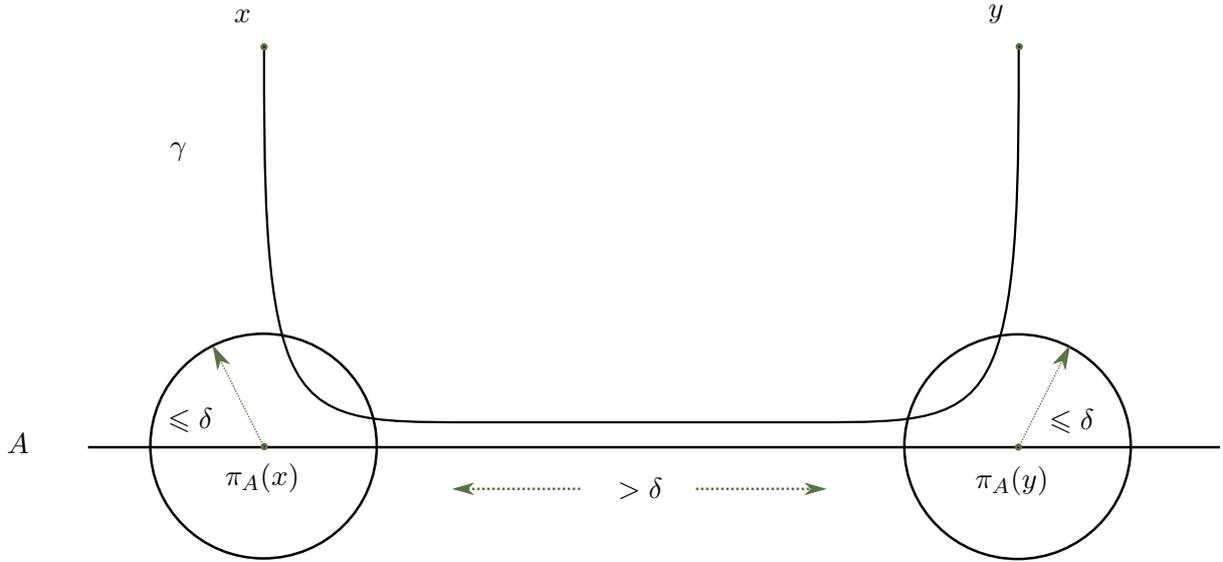
	\caption{The constriction property.}
\end{figure}

\begin{nota}
	Let $\pi_A\colon X\to A$ be a map between $X$ and a subset $A\subset X$. For every $x,y\in X$, we denote $|x-y|_A=|\pi_A(x)-\pi_A(y)|$. For every subset $Y\subset X$, we denote $\diam_A(Y)=\diam(\pi_A(Y))$. For every $x\in X$ and for every pair of subsets $Y,Z\subset X$, we denote $d_A(x,Y)=d(\pi_A(x),\pi_A(Y))$ and $d_A(Y,Z)=d(\pi_A(Y),\pi_A(Z))$. Note that $d_A$ may not be a distance over the collection of subsets of $X$: it may not satisfy the triangle inequality. We will keep this notation for the rest of the paper.
\end{nota}

The following are some standard properties:

\begin{prop}
	\label{prop:elementary properties}
	For every $\delta\ge 0$, there exist a constant $\theta\ge 0$ and a pair of maps, $\sigma \colon \R_{\ge 1}\times \R_{\ge 0}\to \R_{\ge 0}$ and $\zeta \colon \R_{\ge 0}\to \R_{\ge 0}$, such that any $\delta$-constricting map $\pi_A\colon X\to A$ satisfies the following properties:
	
	\begin{enumerate}[label=(\arabic*)]
		\item \textbf{Coarse nearest-point projection.}
		
		\noindent
		For every $x\in X$, we have $|x-\pi_A(x)|\le \mu d(x,A)+\theta.$
		
		\item \textbf{Coarse equivariance.}
		
		\noindent
		Let $H$ be a group acting by isometries on $X$ such that $A$ and $\scrP$ are $H$-invariant. Then for every $h\in H$ and for every $x\in X$, we have $|\pi_A(hx)-h\pi_A(x)|\le \theta$.
		
		\item \textbf{Coarse Lipschitz map.}
		
		\noindent
		For every $x,y\in X$, we have $|x-y|_A\le \mu |x-y|+\theta.$
		
		\item \textbf{Intersection--image.}
		
		\noindent
		For every $\gamma\in \scrP$, we have $| \diam(A^{+\delta}\cap\gamma)-\diam_A(\gamma)|\le \theta.$
		
		\item \textbf{Behrstock inequality.}
		
		\noindent
		Let $\pi_B \colon X\to B$ be a $\delta$-constricting map.  Then for every $x\in X$, we have $$\min\qty{d_A(x,B),d_B(x,A)}\le \theta.$$
		
		\item \textbf{Morseness.}
		
		\noindent
		Let $\kappa\ge 1$, $l\ge 0$. Let $\alpha$ be a $(\kappa,l)$-quasi-geodesic of $X$ with endpoints in $A$. Then $\alpha\subset A^{+\sigma(\kappa,l)}$.
		
		\item \textbf{Coarse invariance.}
		
		\noindent
		Let $\varepsilon\ge 0$. Let $B\subset X$ be a subset such that $d_{\operatorname{Haus}}(A,B)\le \varepsilon$. Then $B$ is $\zeta(\varepsilon)$-constricting.
	\end{enumerate}
\end{prop}

\begin{proof}
	We give some references. For (1), (3) and (4), see \cite[Lemma~2.4]{sisto_contracting_2018}. For (5), see \cite[Lemma~2.5]{sisto_contracting_2018}. For (6), see \cite[Lemma~2.8~(1)]{sisto_contracting_2018}. We leave the proof of the properties (2) and (7) as an exercise. 
\end{proof}

\paragraph{Path system groups.} Let $G$ be a group acting by isometries on a metric space $X$. The \emph{quasi-stabilizer} $\Stab_G(x,r)$ of $x\in X$ of radius $r\ge 0$ is defined as
$$\Stab_G(x,r)=\{g\in G\colon |x-gx|\le r\}.$$ 
The action of $G$ on $X$ is \emph{proper} if for every $x\in X$ and for every $r\ge 0$, we have $|\Stab_G(x,r)|<\infty$. Let $\eta\ge 0$. The action of $G$ on $X$ is $\eta$-\emph{cobounded} if for every $x,x'\in X$, there exists $g\in G$ such that $|x-gx'|\le \eta$.

\begin{df}[Path system group]	
	Let $\mu\ge 1$, $\nu\ge 0$. A $(\mu,\nu)$-\emph{path system group} $(G,X,\scrP)$ is a group $G$ acting properly on a metric space $X$ together with a $G$-invariant collection $\scrP$ of paths of $X$ such that $(X,\scrP)$ is a $(\mu,\nu)$-path system space.
\end{df}

We fix $\mu\ge 1$, $\nu\ge 0$ and a $(\mu,\nu)$-\emph{path system group} $(G,X,\scrP)$.

\begin{df}[Quasi-convex subgroup]
	\label{def:quasi-convexity path system}
	A \emph{subgroup} $H\le G$ is $\eta$-\emph{quasi-convex} if there exists an $H$-invariant $\eta$-quasi-convex subset $Y\subset X$ such that the action of $H$ on $Y$ is $\eta$-cobounded. We will write $(H,Y)$ when we need to stress the $\eta$-quasi-convex subset $Y$ that $H$ is preserving.
\end{df}

\begin{df}[Constricting element]
	\label{def:constricting element path system}
	Let $\delta\ge 0$. An \emph{element} $g\in G$ is $\delta$-\emph{constricting} if the following holds:
	\begin{enumerate}[label=(CE\arabic*)]
		\item $g$ has infinite order.
		\item There exists a $\langle g\rangle$-invariant $\delta$-constricting subset $A\subset X$ so that the action of $\langle g \rangle$ on $A$ is $\delta$-cobounded.
	\end{enumerate}
	We will write $(g,A)$ when we need to stress the $\delta$-constricting subset $A$ that $\langle g \rangle$ is preserving.
\end{df}

\begin{rem}
Note that \autoref{def:quasi-convexity path system} and \autoref{def:constricting element path system} imply the corresponding definitions of the introduction. The converse implication is also true for \autoref{def:constricting element path system}, but the argument requires \autoref{prop:elementary properties} (7) \emph{Coarse invariance}.
\end{rem}

	\section{Growth estimation criteria}
\label{sec:growth-estimation}

In this section, we fix a group $G$ acting properly on a metric space $X$ and a subgroup $H\le G$. The goal is to establish simple criteria so that we can check if $H$ is a solution to the system of equations
\begin{equation*}
	\begin{cases}
		\omega(H)<\omega(G),\\
		\omega(G/H)=\omega(G).
	\end{cases}
\end{equation*}

Our criterion to estimate the relative exponential growth rate is basically \cite[Criterion~2.4]{dalbo_growth_2011}. The statement that we actually need is more specific, so we will give a proof for the convenience of the reader. Recall that the action of a subgroup $H\le G$ on $X$ is \emph{divergent} if its Poincaré series $\mathscr{P}_H(s)$ diverges at $s=\omega(H)$.

\begin{prop}[{\cite[Criterion~2.4]{dalbo_growth_2011}}]
	\label{prop:Dalbo}
	Assume that the following conditions are true:
	\begin{enumerate}[label=(\roman*)]
		\item $\omega(H)<\infty$.
		\item The action of $H$ on $X$ is divergent.
		\item There exist subgroups $K\le G$ and $F\le H\cap K$ so that $F$ is a proper finite subgroup of $K$ and the natural homomorphism  $\phi \colon H\ast_F K\to G$  is injective.
	\end{enumerate}
	Then $\omega(H)<\omega(G)$.
\end{prop}

\begin{rem}
	In the proof below, note that the relative exponential growth rate makes sense for any subset of $G$, as it does the notion of Poincaré series.
\end{rem}

\begin{proof}
	Since the action of $H$ on $X$ is divergent, in particular $H$ is infinite and hence $H-F$ is non-empty. Since $F$ is a proper subgroup of $K$, there exists $k\in K-F$. Denote by $U$ the set of elements of $H\ast_F K$ that can be written as words that alternate elements of $H-F$ and $k$, always with an element of $H-F$ at the beginning and with a $k$ at the end. The inequality  $\omega(\phi(U)) \le \omega(G)$ can be deduced from the definition. It is enough to prove that there exists $s_0\ge 0$ such that $\omega(H)<s_0\le\omega(\phi(U))$. Let $o\in X$. Since $\omega(H)<\infty$, the interval $(\omega(H),\infty)$ is non-empty. Since the action of $H$ on $X$ is divergent, there exists $s_0\in (\omega(H),\infty)$ such that $\sum_{h\in H-F} e^{-s_0 |o-hko|}>1$; otherwise one obtains a contradiction with the divergence of the action of $H$ on $X$.
	
	In order to obtain the inequality $s_0\le \omega(\phi(U))$, it suffices to show that the Poincaré series $\mathscr{P}_{\phi(U)}(s)=\sum_{g\in \phi(U)} e^{-s |o-go|}$ diverges at $s=s_0$. Since $\phi \colon H\ast_F K\to G$ is injective, we have
	$$\mathscr{P}_{\phi(U)}(s)\ge \sum_{m\ge 1}\sum_{h_1,\cdots, h_m\in H-F} e^{-s |o-h_1kh_2k\cdots h_mko|}.$$
	By the triangle inequality, for every $m\ge 1$ and for every $h_1,\cdots, h_m\in H-F$, we have $|o-h_1kh_2k\cdots h_mko|\le \sum_{i=1}^m |o-h_iko|$. Thus,
	$$\sum_{h_1,\cdots, h_m\in H-F} e^{-s |o-h_1kh_2k\cdots h_mko|}\ge \left[\sum_{h\in H-F} e^{-s |o-hko|}\right]^m.$$
	We see that $\mathscr{P}_H(s_0)=\infty$ follows from the claim.
\end{proof}

Our criterion to estimate the quotient exponential growth rate is the following:

\begin{df}
	\label{def:quasi-bijective}
	Let $\phi\colon G\to G$. We say that $G$ is $\phi$-\emph{coarsely} $G/H$ if there exist $\theta\ge 0$ and $x\in X$ satisfying the following conditions:
	\begin{enumerate}[label=(CQ\arabic*)]
		\item For every $u,v\in G$, if $\phi(u)H=\phi(v)H$, then $|\phi(u)x-\phi(v)x|\le \theta$.
		\item For every $u\in G$, $|ux-\phi(u)x|\le \theta$.
	\end{enumerate}
\end{df}

\begin{prop}
	\label{prop:quasi-bijection}
	If there exist $\phi\colon G\to G$ such that $G$ is $\phi$-coarsely $G/H$, then $\omega(G)=\omega(G/H)$.
\end{prop}

\begin{proof}
	The inequality $\omega(G/H)\le \omega(G)$ can be deduced from the defintion. Assume that there exist $\phi\colon G\to G$ such that $G$ is $\phi$-coarsely $G/H$ for $x\in X$ and $\theta\ge 0$.
	\begin{cla}
		There exist $\kappa\ge 1$ such that for every $r>0$, $$|\Stab_G(x,r)|\le \kappa |p(\Stab_G(x,r+\theta))|.$$
	\end{cla}
	Let $\kappa=|\Stab_G(x,3\theta)|$. Let $r>0$. Let $p\colon G \onto G/H$ be the natural projection. Let $q\colon G \rightarrow G/H$ the map that sends $u$ to $\phi(u)H$. Note that the quasi-stabilizer $\Stab_G(x,r)$ can be decomposed as the disjoint union of the sets $q^{-1}(q(u))$ such that $q(u) \in q(\Stab_G(x,r))$. Hence,
	$$|\Stab_G(x,r)|\le \sum_{q(u)\in q(\Stab_G(x,r))}| q^{-1}(q(u))|.$$
	It suffices to estimate the size of $q(\Stab_G(x,r))$ and the size of $q^{-1}(q(u))$, for every $u\in G$. First we prove that $|q(\Stab_G(x,r))|\le |p(\Stab_G(x,r+\theta))|$. Let $u\in \Stab_G(x,r)$. By the triangle inequality,
	$$|x-\phi(u)x|\le |x-ux|+|ux-\phi(u)x|.$$
	By the hypothesis (CQ2), we have $|ux-\phi(u)x|\le \theta$. Hence $|x-\phi(u)x|\le r+\theta$. Consequently, $q(\Stab_G(x,r))\subset p(\Stab_G(x,r+\theta))$. Now we prove that for every $u\in G$, we have $|q^{-1}(q(u))| \le \kappa$. Let $u\in G$. Since $|u\Stab_G(x,3\theta)|=|\Stab_G(x,3\theta)|=\kappa$, it is enough to prove that $u^{-1}q^{-1}(q(u))\subset \Stab_G(x,3\theta)$. Let $v\in q^{-1}(q(u))$. 
	By the triangle inequality,
	\begin{equation*}
		\begin{aligned}
			|x-u^{-1}vx|=|ux-vx|\le |ux-\phi(u)x|+|\phi(u)x-\phi(v)x|+|\phi(v)x-vx|.
		\end{aligned}
	\end{equation*}
	Since $q(u)=q(v)$, we have that $\phi(u)H=\phi(v)H$. It follows from the hypothesis (CQ1) that $|\phi(u)x-\phi(v)x|\le \theta$. By the hypothesis (CQ2), we have $\max\{|ux-\phi(u)x|,|vx-\phi(v)x|\}\le\theta$. Thus, $|x-u^{-1}vx|\le 3\theta$. This proves the claim. 
	
	Consequently,
	$$\omega(G)\le \limsup_{r\to \infty}\frac{1}{r} \log |p(\Stab_G(x,r+\theta))|.$$
	Finally, observe that
	$$\limsup_{r\to \infty}\frac{1}{r} \log |p(\Stab_G(x,r+\theta))|= \limsup_{r\to \infty}\frac{r+\theta}{r}\frac{1}{r+\theta} \log |p(\Stab_G(x,r+\theta))|.$$
	Hence $\omega(G)\le \omega(G/H)$.
\end{proof}

	\section{Buffering sequences}
\label{sec:buffering}
In this section, we fix constants $\mu\ge 1$, $\nu\ge 0$ and a $(\mu,\nu)$-path system space $(X,\scrP)$. Despite the fact that our space $X$ does not carry any global geometric condition, we still can obtain some control through constricting subsets. We could ignore the ``wild regions'' if, for instance, we were able to ``jump'' from one constricting subset to another. The buffering sequences  below encapsulate this idea. In fact, the proofs of our main results consist essentially in building up some particular buffering sequences. W. Yang had already introduced this concept for piece-wise geodesics in \cite{yang_statistically_2019}.

\begin{df}
	Let $\delta$, $\varepsilon$, $L\ge 0$. Let $\mathscr{A}$ be a collection of subsets of $X$. A finite sequence of subsets $Y_0, A_1, Y_1,\cdots, A_n, Y_n\subset X$ where $Y_0$ and $Y_n$ are the only possible empty sets is $(\delta,\varepsilon,L)$-\emph{buffering on} $\mathscr{A}$ if for every $i\in\zinterval{1}{n}$ the set $A_i$ belongs to $\mathscr{A}$ and there exists a $\delta$-constricting map $\pi_{A_i}\colon X\to A_i$ with the following properties whenever $Y_i$ and $Y_{i-1}$ are non-empty:
	\begin{enumerate}[label=(BS\arabic*)]
		\item $\max\{\diam_{A_i}(A_{i+1}),\diam_{A_{i+1}}(A_{i})\}\le \varepsilon$ if $i\neq n$.
		\item $\max\{\diam_{A_i}(Y_{i-1}),\diam_{A_i}(Y_{i})\}\le \varepsilon$.
		\item $\max\{d(A_i,Y_{i-1}),d(A_i,Y_{i})\}\le \varepsilon$.
		\item $d_{A_i}(Y_{i-1},Y_i)\ge L$.
	\end{enumerate}
\end{df}

\begin{figure}[htbp]
	\centering
	\includegraphics[width=0.6\textwidth]{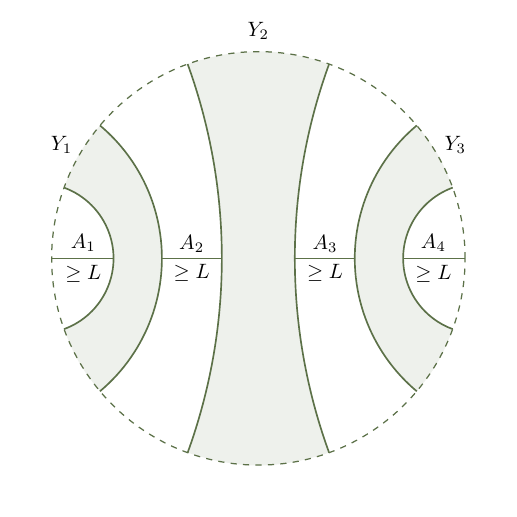}
	\caption{An example of a 
		buffering sequence in the Poincaré disk model. In this example, the sets $A_i$ are subpaths of length $\ge L$ of a given bi-infinite geodesic $\alpha$. Each set $Y_i$ is the collection of geodesics that are orthogonal to the geodesic segment of $\alpha$ that is between $A_i$ and $A_{i+1}$. In particular, the sets $Y_i$ are quasi-convex. For more intuition, one could interpret this picture on a tree.}
\end{figure}

What makes buffering sequences remarkable is that they satisfy a variant of \emph{Behrstock inequality}. We will find a direct application of the following inequality later in the study of the quotient exponential growth rates:

\begin{prop}
	\label{prop:Behrstock's inequality 2}
	For every $\delta$, $\varepsilon\ge 0$, there exists $\theta\ge 0$ with the following property. Let $A,Y,B\subset X$ be a $(\delta,\varepsilon,0)$-buffering sequence on $\{A,B\}$. Then for every $x\in X$,
	$$\min\qty{d_{A}(x,Y),d_{B}(x,Y)}\le \theta.$$
\end{prop}

\begin{proof}
	Let $\delta$, $\varepsilon\ge 0$. Let $\theta_0=\theta_0(\delta)\ge 0$ be the constant of \autoref{prop:elementary properties}. Let $\theta> \theta_0+1$. Its exact value will be precised below. Let $A,Y,B\subset X$ be a $(\delta,\varepsilon,0)$-buffering sequence on $\{A,B\}$. Let $x\in X$. By symmetry, it suffices to show that if $d_A(x,Y)>\theta$, then $d_B(x,Y)\le \theta$. Assume that $d_A(x,Y)>\theta$. Let $a\in A$ such that $|x-a|_B\le d_B(x,A)+1$. Let $b\in B$. Let $y\in Y$. By (BS3), we have $\max\{d(A,Y),d(B,Y)\}\le \varepsilon$; hence there exist $p \in A^{+\varepsilon+1}\cap Y$ and $q \in B^{+\varepsilon+1}\cap Y$.  It follows from the definition of buffering sequence that
	$$\max\qty{|b-\pi_{B}(q)|_A, |q-p|_A,|a-\pi_{A}(p)|_B, |p-y|_B} \le \varepsilon.$$
	Applying together \autoref{prop:elementary properties} (1) \emph{Coarse nearest-point projection} and (3) \emph{Coarse Lipschitz map}, we obtain
	$$\max\{|\pi_{B}(q)-q|_A,|\pi_{A}(p)-p|_B\}\le \mu^2(\varepsilon+1)+\mu \theta_0+\theta_0.$$
	\begin{cla}
		$d_{A}(x,B)>\theta_0$
	\end{cla}
	By the triangle inequality,
	$$|x-b|_A\ge |x-p|_A-|b-\pi_{B}(q)|_A-|\pi_{B}(q)-q|_A-|q-p|_A.$$
	Moreover, $|x-p|_A\ge d_A(x,Y)$. Since the element $b$ is arbitrary and we have $d_A(x,Y)> \theta_0+1$, we obtain $d_{A}(x,B)>\theta_0$. This proves the claim.
	
	Finally, we are going to estimate $d_B(x,Y)$. By the triangle inequality,
	$$|x-y|_B\le |x-a|_B+|a-\pi_{A}(p)|_B+|\pi_{A}(p)-p|_B+|p-y|_B.$$
	Since $d_A(x,B)>\theta_0$, it follows from \autoref{prop:elementary properties} (5) \emph{Behrstock inequality} and the definition of $a$ that $|x-a|_B\le \theta_0+1$. Since the element $y$ is arbitrary, we obtain $d_{B}(x,Y)\le \theta$ for $\theta=2\theta_0+1+2\varepsilon+\mu^2(\varepsilon+1)+\mu\theta_0$.
\end{proof}

The corollary below will be applied to the study of the relative exponential growth rates:
\begin{cor}
	\label{cor:buffer sequence}
	For every $\delta$, $\varepsilon$, $\theta\ge 0$ there exists $L\ge 0$ with the following property. Let $Y_0, A_1, Y_1,\cdots, A_n, Y_n\subset X$ be an $(\delta,\varepsilon,L)$-buffering sequence on $\{A_i\}$. Then for every $i\in\zinterval{1}{n}$, 
	$$d_{A_{i}}(Y_0,Y_{i})>\theta.$$
\end{cor}

\begin{proof}
	Let $\delta$, $\varepsilon$, $\theta\ge 0$. Let $\theta_0=\theta_0(\delta,\varepsilon)\ge 0$ be the constant of \autoref{prop:Behrstock's inequality 2}. We put $L=\theta+\theta_0+1$. Let $y_0\in Y_0$. Let $i\in\zinterval{1}{n}$. 
	\begin{cla}
		$d_{A_{i}}(y_0,Y_{i})\ge d_{A_{i}}(Y_{i-1},Y_{i})-d_{A_{i}}(y_0,Y_{i-1})$.
	\end{cla}
	Let $y_{i-1}\in Y_{i-1}$ and $y_{i}\in Y_{i}$. By the triangle inequality,
	$$|y_0-y_{i}|_{A_i}\ge |y_{i-1}-y_{i}|_{A_i}-|y_0-y_{i-1}|_{A_i}.$$
	Note that $|y_{i-1}-y_{i}|_{A_i}\ge d_{A_{i}}(Y_{i-1},Y_{i})$. Since the elements $y_{i-1}, y_i$ are arbitrary, this proves the claim.
	
	Finally, we prove by induction on $i\in\zinterval{1}{n}$ that, $d_{A_{i}}(Y_0,Y_{i})>\theta$. If $i=1$, then $d_{A_{1}}(Y_0,Y_{1})>\theta$ follows from (BS4), since $L>\theta$. Assume that $i\in\zinterval{1}{n-1}$ and $d_{A_{i}}(Y_0,Y_{i})>\theta$. Then $d_{A_{i}}(y_0,Y_{i})>\theta_0$. It follows from \autoref{prop:Behrstock's inequality 2} that $d_{A_{i+1}}(y_0,Y_{i})\le \theta_0$. By (BS4), $d_{A_{i+1}}(Y_{i},Y_{i+1})\ge L$. Applying the previous claim, we obtain $d_{A_{i+1}}(y_0,Y_{i+1})> \theta$. Since the element $y_0$ is arbitrary, $d_{A_{i+1}}(Y_0,Y_{i+1})> \theta$. This concludes the inductive step.
\end{proof}

	\section{Quasi-convexity in the intersection--image property}
\label{sec:qc-int-img}

In this section, we fix constants $\mu\ge 1$, $\nu\ge 0$ and a $(\mu,\nu)$-path system space $(X,\scrP)$. In this section, we prove a variant of \autoref{prop:elementary properties} (4) \emph{Intersection--Image}. Basically, we will be exchanging paths of $\scrP$ for quasi-convex subsets of $X$, further thickening the involved sets.

\begin{prop}
	\label{prop:diameter}
	For every $\delta$, $\eta\ge 0$, there exist $\theta\ge 0$ and $\zeta\colon \R_{\ge0}\times \R_{\ge 0}\to\R_{\ge 0}$ with the following property. Let $\pi_A\colon X\to A$ be a $\delta$-constricting map. Let $Y$ be an $\eta$-quasi-convex subset of $X$. Let $\varepsilon_1\ge 0$, $\varepsilon_2\ge 0$. Then 
	$$|\diam(A^{+\theta+\varepsilon_1}\cap Y^{+\varepsilon_2})-\diam_A(Y)|\le \zeta(\varepsilon_1,\varepsilon_2).$$
\end{prop}

\begin{proof}
	Let $\delta$, $\eta\ge 0$. Let $\theta_0=\theta_0(\delta)\ge 0$ be the constant of \autoref{prop:elementary properties}. We put $\theta=\delta+\eta+1$. Let $\zeta\colon \R_{\ge0}\times \R_{\ge 0}\to\R_{\ge 0}$ depending on $\delta,\eta$. Its exact value will be precised below. Let $\pi_A\colon X\to A$ be a $\delta$-constricting map. Let $Y$ be an $\eta$-quasi-convex subset of $X$. Let $\varepsilon_1\ge 0$, $\varepsilon_2\ge 0$.
	
	First we prove that $\diam_A(Y)\le \diam(A^{+\theta+\varepsilon_1}\cap Y^{+\varepsilon_2})+\zeta(\varepsilon_1,\varepsilon_2)$. Let $x,y\in Y$. It suffices to assume that $|x-y|_A>\delta$. Let $\gamma\in\scrP$ joining $x$ to $y$. By (CS2), there exist $p,q\in\gamma$ such that
	$$\max\{|\pi_A(x)-p|,|\pi_A(y)-q|\}\le \delta.$$
	Since the subset $Y$ is $\eta$-quasi-convex, there exist $p',q'\in Y$ such that
	$$\max\{|p-p'|,|q-q'|\}\le \eta+1.$$
	By the triangle inequality,
	$$|x-y|_A\le |\pi_A(x)-p|+|p-p'|+|p'-q'|+|q'-q|+|q-\pi_A(y)|.$$
	Since $p',q'\in A^{+\theta+\varepsilon_1}\cap Y^{+\varepsilon_2}$, we have $|p'-q'|\le \diam(A^{+\theta+\varepsilon_1}\cap Y^{+\varepsilon_2})$. Hence,
	$$|x-y|_A\le \diam(A^{+\theta+\varepsilon_1}\cap Y^{+\varepsilon_2})+2\delta+2\eta+1.$$
	
	Now we prove that $\diam(A^{+\theta+\varepsilon_1}\cap Y^{+\varepsilon_2})\le \diam_A(Y)+\zeta(\varepsilon_1,\varepsilon_2)$. Let $x,y\in A^{+\theta+\varepsilon_1}\cap Y^{+\varepsilon_2}$. Since $x,y\in Y^{+\varepsilon_2}$, there exist $x',y'\in Y$ such that $\max\{|x-x'|,|y-y'|\}\le \varepsilon_2+1$. By the triangle inequality,
	$$|x-y|\le |x-\pi_A(x)|+|x-x'|+|x'-y'|_A+|y'-y|_A+|\pi_A(y)-y|.$$
	Since $x,y\in A^{+\theta+\varepsilon_1}$, it follows from \autoref{prop:elementary properties} (1) \emph{Coarse nearest-point projection} that
	$$\max\{|x-\pi_A(x)|,|y-\pi_A(y)|\}\le \mu(\theta+\varepsilon_1)+\theta_0.$$
	It follows from \autoref{prop:elementary properties} (3) \emph{Coarse Lipschitz Map} that,
	$$\max\{|x-x'|_A,|y-y'|_A\}\le \mu(\varepsilon_2+1)+\theta_0.$$
	Since $\pi_A(x'),\pi_A(y')\in \pi_A(Y)$, we have $|x'-y'|_A\le \diam_A(Y)$. Hence,
	$$|x-y|\le \diam_A(Y)+2\mu(\theta+\varepsilon_1)+2\mu(\varepsilon_2+1)+4\theta_0.$$
	Finally, we put $\zeta(\varepsilon_1,\varepsilon_2)=\max\{2\delta+2\eta+1,2\mu(\theta+\varepsilon_1)+2\mu(\varepsilon_2+1)+4\theta_0\}$.
\end{proof}
Applying the symmetry of \autoref{prop:diameter} in combination with \autoref{prop:elementary properties} (6) \emph{Morseness} and  (7) \emph{Coarse invariance}, we deduce:
\begin{cor}
	\label{cor:symmetry in projections}
	For every $\delta\ge 0$, there exists $\theta\ge 0$ with the following property. Let $\pi_A\colon X\to A$ and $\pi_B\colon X\to B$ be $\delta$-constricting maps. Then:
	$$|\diam_{A}(B)-\diam_{B}(A)|\le \theta.$$
\end{cor}

	\section{Finding a quasi-convex element}
\label{sec:finding}
Given a torsion-free hyperbolic group $G$ containing a loxodromic element $g_0$ and an infinite index quasi-convex subgroup $H$, one can find another loxodromic element $g\in G$ conjugate to $g_0$ so that $H$ has trivial intersection with $\langle g\rangle$ \cite[Theorem 1]{arzhantseva_quasiconvex_2001}. The goal of this section is to reimplement this fact in our setting, using a ``quasi-convex element'' instead of a loxodromic element. 

\begin{conv}
In this section, we fix:
\begin{itemize}[label=$\blacktriangleright$]
	\item Constants $\mu\ge 1$, $\nu\ge 0$.
	\item A $(\mu,\nu)$-path system group $(G,X,\scrP)$.
\end{itemize}
\end{conv}

\begin{df}[Quasi-convex element]
	Let $\eta\ge 0$. An \emph{element} $g\in G$ is $\eta$-\emph{quasi-convex} if the following holds:
	\begin{enumerate}[label=(QE\arabic*)]
		\item $g$ has infinite order.
		\item $\langle g \rangle$ is an $\eta$-quasi-convex subgroup of $G$.
	\end{enumerate}
	We will write $(g,A)$ when we need to stress the $\eta$-quasi-convex subset $A$ that $\langle g \rangle$ is preserving.
\end{df}

The main result of this section is the following.

\begin{prop}
	\label{prop:alternative of intersection and hausdorff distance}
	Let $\eta\ge 0$. Assume that $G$ contains an $\eta$-quasi-convex element $(g,A)$. There exists $\theta=\theta(\eta,g,A) \ge 1$ satisfying the following. Let $(H,Y)$ be an $\eta$-quasi-convex subgroup of $G$. Then:
	\begin{enumerate}[label=(\roman*)]
		\item For every $u\in G$, if $\diam(uA\cap Y)> \theta$, then $uA\subset Y^{+\theta}$.
		\item Let $H\le K\le G$. If $[K:H]>\theta$, then there exist $k\in K$ such that $\diam(kA\cap Y)\le \theta$.
	\end{enumerate}
\end{prop}
\begin{rem}
	Under the notation of (ii), when $K=G$, the element $kgk^{-1}$ has the desired property that we were looking for. Note that $(kgk^{-1},kA)$ is quasi-convex since $\scrP$ is $G$-invariant.
\end{rem}

The rest of the section is devoted to the proof of \autoref{prop:alternative of intersection and hausdorff distance}.

\begin{df}
	Let $\kappa\ge 1$, $l\ge 0$. A map $\phi\colon (Y,d_Y)\to (Z,d_Z)$ between two metric spaces is a $(\kappa,l)$-\emph{quasi-isometric embedding} if for every $y,y'\in Y$, $$\frac{1}{\kappa} d_Y(y,y')-l\le d_Z(\phi(y),\phi(y'))\le \kappa d_Y(y,y')+l.$$
\end{df}

We start with a variant of Milnor-Schwarz Theorem. If $U$ is a generating set of a group $H$, we denote by $d_U$ the word metric of $H$ with respect to $U$.

\begin{lem}
	\label{lem:Milnor-Schwarz}
	For every $\eta\ge 0$, there exist $\theta \ge 1$ with the following property. Let $(H,Y)$ be an $\eta$-quasi-convex subgroup of $G$. For every $y\in Y$, there exists a finite generating set $U$ of $H$ such that the orbit map $(H,d_{U})\to X$, $h\mapsto hy$ is a $(\theta,\theta)$-quasi-isometric embedding.
\end{lem}

For the proof, one can use the same kind of argument as that of Milnor-Schwarz Theorem, but bearing in mind that $Y$ might not be a length metric space, which is required by the original statement. The only difference here is that one uses the paths of $\scrP$ with endpoints in $Y$. They are enough for the proof since they approximate sufficiently well the distances, at least in this situation.

\begin{lem}
	\label{lem:abelian}
	Let $\eta\ge 0$. Let $H\le G$ be an abelian subgroup. Let $Y\subset X$ be an $H$-invariant subset so that the action of $H$ on $Y$ is $\eta$-cobounded. Then, for every $h\in H$ and for every $y,z\in Y$,
	$$\big||y-hy|-|z-hz|\big|\le 2\eta.$$
\end{lem}

\begin{proof}
	Let $h\in H$. Let $y,z\in Y$. Since the action of $H$ on $Y$ is $\eta$-cobounded, there exists $k\in H$ such that $|z-ky| \le \eta$. By the triangle inequality,
	$$|y-hy|\le |ky-khy|\le |ky-z|+|z-hz|+|hz-khy|.$$
	Since the subgroup $H$ is abelian, $|hz-khy|=|z-ky|$.
	Thus, $|y-hy|\le |z-hz|+2\eta$. Finally, exchanging the roles of $y$ and $z$, we obtain $|y-hy|\ge |z-hz|-2\eta$. 
\end{proof}

Next, we are going to check that we can obtain uniform quasi-isometric embeddings of $\Z$ in $X$ via the orbit maps of quasi-convex elements of $G$ that share the same constant. For this reason, we introduce the following definition:

\begin{df}
	Let $g\in G$. Let $x\in X$. The \emph{stable translation length} of $g$ is 
	$$\stlen{g}=\limsup_{m\to\infty}\frac{1}{m}\, |g^mx-x|.$$
	Note that $\stlen{g}$ does not depend on the choice of the point $x\in X$.
\end{df}

\begin{rem}
	\label{rk:asymptotic translation length}
	Let $g\in G$. By subadditivity, for every $x\in X$, we have 
	$$\stlen{g}=\inf_{m\ge 1} \frac{1}{m} |g^mx-x|=\lim_{m\to\infty}\frac{1}{m}\, |g^mx-x|.$$
\end{rem}

\begin{lem}
	\label{lem:characterisation of qi embedding}
	Let $\eta\ge 0$. Let $g\in G$. Let $A\subset X$ be a $\langle g\rangle$-invariant subset so that the action of $\langle g \rangle$ on $A$ is $\eta$-cobounded. The following statements are equivalent:
	\begin{enumerate}[label=(\roman*)]
		\item There exists $x\in X$ such that the orbit map $\Z\to X$, $m\mapsto g^mx$ is a quasi-isometric embedding.
		\item  $\stlen{g}>0$.
		\item There exists $\theta=\theta(\eta,g,A)\ge 1$ such that for every $a\in A$, the orbit map $\Z\to X$, $m\mapsto g^ma$ is a $(\theta,0)$-quasi-isometric embedding.
	\end{enumerate}
\end{lem}

\begin{proof}
	The implication $(iii) \Rightarrow (i)$ already holds.
	
	$(i) \Rightarrow (ii)$. Assume that there exists $x\in X$ such that the orbit map $\Z\to X$, $m\mapsto g^mx$ is a quasi-isometric embedding. Then there exist $\kappa\ge 1$, $l \ge 0$ such that for every $m\ge 1$,
	$$\frac{1}{\kappa}-\frac{l}{m}\le \frac{1}{m}|x-g^mx|\le \kappa+\frac{l}{m}.$$
	Therefore, $\stlen{g} \ge \frac{1}{\kappa}>0$.
	
	$(ii) \Rightarrow (iii)$. Assume that $\stlen{g}>0$. Let $\|g\|_A=\inf_{a\in A}|a-ga|$. Then we can define $\theta=\max\left\{\|g\|_A+2\eta,\frac{1}{\stlen{g}},1\right\}$. Let $a\in A$. Applying the triangle inequality we obtain that for every $m\in \Z$, $|a-g^ma|\le |a-ga||m|$. It follows from \autoref{lem:abelian} that $|a-ga|\le \|g\|_A+2\eta$. Since $\stlen{g}=\inf_{n\in\Z -\{0\}} \frac{1}{|n|}|a-g^{|n|}a|$, we obtain that for every $m\in \Z$, $|a-g^ma|\ge \stlen{g}|m|$. Hence the orbit map $\Z\to X$, $m\mapsto g^ma$ is a $(\theta,0)$-quasi-isometric embedding.
\end{proof}

\begin{lem}
	\label{lem:qi embedding of quasiconvex subgroups}
	Let $\eta\ge 0$. Let $(g,A)$ be an $\eta$-quasi-convex element of $G$. There exists $\theta=\theta(\eta,g,A)\ge 1$ such that for every $a\in A$, the orbit map $\Z\to X$, $m\mapsto g^ma$ is a $(\theta,0)$-quasi-isometric embedding. Moreover, $\stlen{g}>0$.
\end{lem}

\begin{proof}
	We are going to apply \autoref{lem:Milnor-Schwarz} and \autoref{lem:characterisation of qi embedding}. Let $a\in A$. According to \autoref{lem:Milnor-Schwarz}, there exist a finite generating set $U$ of $\langle g\rangle $ such that the orbit map $\phi \colon (\langle g \rangle,d_{U})\to X$, $h\mapsto ha$ is a quasi-isometric embedding. Furthermore, since $g$ has infinite order, the map $\chi \colon \Z \to  \langle g\rangle$, $m\mapsto g^m$ is an isomorphism. Let $V=\chi^{-1}(U)$. In particular $\chi \colon (\Z, d_V) \to  (\langle g\rangle, d_U)$  is an isometry. Morover, the map $\psi\colon \Z	\to (\Z, d_V)$ is a quasi-isometric embedding. Hence, the composition $\phi\circ \chi\circ \psi$ is a quasi-isometric embedding. Now both of the statements of the lemma follow from \autoref{lem:characterisation of qi embedding}.
\end{proof} 

We continue by upper bounding the length of a quasi-geodesic of $X$ by the number of points of an orbit of a subgroup $H$ of $G$ that fall inside a precise neighbourhood of this quasi-geodesic, whenever the quasi-geodesic falls also inside a neighbourhood of that orbit.

\begin{lem}
	\label{lem:number of orbit points are minored by length}
	For every $\eta\ge 0$, $\kappa\ge 1$, $l\ge 0$, there exists $\theta\ge 1$ with the following property.  Let $H\le G$. Let $Y\subset X$ be an $H$-invariant subset such that the action of $H$ on $Y$ is $\eta$-cobounded.  Let $y\in Y$. Let $\gamma$ be a $(\kappa,l)$-quasi-geodesic of $X$ such that $\gamma\subset Y^{+\eta}$. Let $U=\{u\in H \colon \, uy\in \gamma^{+2\eta+1}\}$. Then 
	$$\ell(\gamma)\le \theta |U|.$$
\end{lem}

\begin{proof}
	Let $\eta\ge 0$, $\kappa\ge 1$, $l\ge 0$. Let $\theta= \theta(\eta,\kappa,l)\ge 1$. Its exact value will be precised below. Let $H$, $Y$, $y$, $\gamma\colon[0,L]\to X$ and $U$ as in the statement. Let $m=\left\lfloor \frac{L}{\theta} \right\rfloor+1$. We fix a partition $0=t_0\le t_1\le \cdots \le t_m=L$ of $[0,L]$ such that $|t_{m-1}-t_m|\le \theta$ and such that if $m\ge 2$, then for every $i\in\zinterval{0}{m-2}$, we have $|t_i-t_{i+1}|=\theta$. Hence $\ell(\gamma)=L\le \theta m$. We prove that $m\le |U|$. Let $i\in\zinterval{0}{m-1}$. Denote $x_i=\gamma(t_i)$. Since the action of $H$ on $Y$ is $\eta$-cobounded and $\gamma\subset Y^{+\eta}$, for every $i\in\zinterval{0}{m-1}$, there exists $h_i\in H$ such that $|x_i-h_iy|\le 2\eta+1$. In particular, $h_i\in U$. From now on we may assume that $m\ge 2$, otherwise there is nothing to show. Let $i,j\in\zinterval{0}{m-1}$ such that $i\neq j$. We claim that $h_i\neq h_j$. The claim will follow when we show that $|h_iy-h_{j}y|>0$. By the triangle inequality,
	$$|h_iy-h_{j}y|\ge |x_i-x_j|-|x_i-h_iy|-|x_j-h_jy|.$$
	Since $\gamma$ is a $(\kappa,l)$-quasi-geodesic,
	$$|x_i-x_j|\ge \frac{1}{\kappa}|t_i-t_{j}|-\frac{l}{\kappa}.$$
	Since $i,j\in\zinterval{0}{m-1}$, we have that $|t_i-t_{j}|\ge \theta$. To sum up,
	$$|h_iy-h_jy|\ge \frac{\theta}{\kappa}-\frac{l}{\kappa}-4\eta-2.$$
	Finally, we put $\theta=\kappa\left(\frac{l}{\kappa}+4\eta+2\right)+1$. Hence, $|h_iy-h_{j}y|>0$. In particular, we obtain $m\le |U|$.
\end{proof}

The following fact is a direct consequence of the triangle inequality:
\begin{lem}
	\label{lem:conjugation of balls}
	Let $\eta\ge 0$. Let $H\le G$. Let $Y\subset X$ be an $H$-invariant subset so that the action of $H$ on $Y$ is $\eta$-cobounded. Then, for every $y,z\in Y$, there exists $h\in H$ such that for every $r>0$, 
	$$h^{-1}\Stab_G(y,r)h\subset \Stab_G(z,r+2\eta).$$
\end{lem}

Finally, we show that there is a uniform threshold that ensures the existence of a uniformly short element in the intersection of any pair of quasi-convex subgroups of $G$ that share the same constant.

\begin{lem}
	\label{lem:element and intersection}
	For every $\eta\ge 0$, there exists $\theta\ge 1$ with the following property. Let $(H,Y)$ and $(K,Z)$ be $\eta$-quasi-convex subgroups of $G$. If $\diam(Y\cap Z)>\theta$, then there exist $y\in Y\cap Z$ and $h\in H\cap K\cap \Stab_G(y,\theta)-\{1_G\}$.
\end{lem}

\begin{proof}
	Let $\eta\ge 0$. Let $\theta_0=\theta_0(\eta,\mu,\nu)\ge 1$ be the constant of \autoref{lem:number of orbit points are minored by length}. Let $o\in Y$. We denote $W=\Stab_G(o,6\eta+2)$. Let $\theta_1=\theta_0|W|+\theta_0$. Note that the constant $\theta_1$ is finite since the action of $G$ on $X$ is proper. We put $\theta=2\theta_1+4\eta+2$. Let $(H,Y)$ and $(K,Z)$ be $\eta$-quasi-convex subgroups of $G$. Assume that $\diam(Y\cap Z)>\theta$. Since $\diam(Y\cap Z)>\theta_1$, there exist $y,z\in Y\cap Z$ such that $|y-z|>\theta_1$. Let $\beta\in \scrP$ joining $y$ to $z$. Since $\ell(\beta)>\theta_1$, there exist $z'\in \beta$ and a subpath $\gamma$ of $\beta$ joining $y$ to $z'$ such that $\ell(\gamma)=\theta_1$. We denote $U=\{u\in H\colon\, uy\in  \gamma^{+2\eta+1}\}$ and $V=\Stab_G(y,4\eta+2)$. 
	
	The first step is to construct a map $\phi \colon U \to V$. Let $u\in U$. By definition of $U$, there exists $x\in \gamma$ such that $|uy-x|\le 2\eta+1$. Since the subgroup $(K,Z)$ is $\eta$-quasi-convex, there exists $k_u\in K$ such that $|x-k_uy|\le 2\eta+1$. By the triangle inequality,
	$$|uy-k_uy|\le |uy-x|+|x-k_uy|.$$ 
	Consequently, $|u^{-1}k_uy-y|\le 4\eta+2$. Hence, $u^{-1}k_u\in V$. We define $\phi \colon U \to V$ to be the map that sends every $u\in U$ to $u^{-1}k_u\in V$. 
	
	Next, we show that the map $\phi \colon U \to V$ is not injective. Since $Y$ is $\eta$-quasi-convex, we have that $\gamma\subset \beta\subset Y^{+\eta}$. It follows from \autoref{lem:number of orbit points are minored by length} that $|U|\ge \frac{1}{\theta_0}\ell(\gamma)$. By hypothesis, $\ell(\gamma)=\theta_0|W|+\theta_0$. Since the action of $H$ on $Y$ is $\eta$-cobounded, it follows from \autoref{lem:conjugation of balls} that there exists $h\in H$ such that $h^{-1}Vh\subset W$ and hence $|W|\ge|h^{-1}Vh|= |V|$. Consequently, $|U|>|V|$. Therefore, the map $\phi \colon U \to V$ is not injective.
	
	Now we claim that $U\subset \Stab_G(y,\theta_1+2\eta+1)$. Let $u\in U$. By definition of $U$, there exists $x\in\gamma$ such that $d|x-uy|\le 2\eta+1$. By the triangle inequality,
	$$|y-uy|\le |y-x|+|x-uy|.$$
	Moreover, $|y-x|\le \ell(\gamma)= \theta_1$. Hence $|y-uy|\le \theta_1+2\eta+1$.
	
	Finally, since the map $\phi \colon U \to V$ is not injective, there exist $u_1,u_2\in U$ such that $u_1\neq u_2$ and $u_1^{-1}k_{u_1}=u_2^{-1}k_{u_2}$. In particular, $u_2 u_1^{-1}\in H\cap K-\{1_G\}$. Further, according to the triangle inequality,
	$$|y-u_2u_1^{-1}y|\le |y-u_2y|+|u_2y-u_2u_1^{-1}y|.$$
	It follows from the claim above that $|y-u_2u_1^{-1}y|\le \theta$. Therefore, $u_2u_1^{-1}\in H\cap K\cap \Stab_G(y,\theta)-\{1_G\}$.
\end{proof}

We are ready to prove the proposition:

\begin{proof}[ Proof of \autoref{prop:alternative of intersection and hausdorff distance}]
	
	Let $\eta\ge 0$. Assume that $G$ contains an $\eta$-quasi-convex element $(g,A)$. We are going to determine the value of $\theta=\theta(\eta,g,A)\ge 1$. By \autoref{lem:qi embedding of quasiconvex subgroups}, there exists $\theta_0=\theta_0(\eta,g,A)\ge 1$ such that for every $a\in A$, the orbit map $\Z\to X$, $m\mapsto g^ma$ is a $(\theta_0,0)$-quasi-isometric embedding. Let $\theta_1=\theta_1(\eta)\ge 1$ be the constant of \autoref{lem:element and intersection}. Let $\theta_2=\eta+\theta_0^2\theta_1$.  Let $o\in A$. We denote $U=\Stab_G(o,2\theta_2+\eta+1)$. Let $\theta=\max\{\theta_2, |U|\}$. Note that the constant $\theta$ is finite since the action of $G$ on $X$ is proper. Let $(H,Y)$ be an $\eta$-quasi-convex subgroup of $G$.
	
	\begin{enumerate}[label=(\roman*)]
		\item Let $u\in G$. Assume that $\diam(uA\cap Y)>\theta$. Let $a\in A$. We prove that $ua \in Y^{+\theta_2}$. Since $\scrP$ is $G$-invariant, the element $(ugu^{-1},uA)$ is $\eta$-quasi-convex. Since $\diam(uA\cap Y)>\theta_1$, according to \autoref{lem:element and intersection}, there exist $b\in A$ and $M\in \Z-\{0\}$ such that $ub\in uA\cap Y$ and $ug^{M}u^{-1}\in H\cap \Stab_G(ub,\theta_1)$. Since the action of $\langle g \rangle $ on $A$ is $\eta$-cobounded, there exists $m\in \Z$ such that $|a-g^mb|\le \eta$. By Euclid's division Lemma, there exist $q,r\in \Z$ such that $m=qM+r$ and $0\le r\le |M|-1$. By the triangle inequality,
		$$d(ua,Y)\le |ua-ug^{qM}b| \le |ua-ug^mb|+|ug^mb-ug^{qM}b|.$$
		Note that $|ua-ug^mb|=|a-g^mb|\le \eta$. Moreover, it follows from \autoref{lem:qi embedding of quasiconvex subgroups} that 
		$$|ug^mb-ug^{qM}b|=|g^{r}b-b|\le \theta_0|r|.$$
		Note also that $|r|\le |M|$. Applying again \autoref{lem:qi embedding of quasiconvex subgroups}, we obtain that $|M|\le \theta_0 |g^Mb-b|$. By \autoref{lem:element and intersection}, $|g^Mb-b|=|ug^Mu^{-1}ub-ub|\le \theta_1$. Hence,
		$$d(ua,Y)\le \theta_2\le \theta.$$
		
		\item  Let $H\le K\le G$. We argue by contraposition. Assume that for every $k\in K$, we have $\diam(kA\cap Y)>\theta$. We prove that $[K:H]\le |U|$. It follows from (i) that $KA\subset Y^{+\theta_2}$. Then there exists $y\in Y$ such that $|o-y|\le \theta_2+1$. Since the action of $H$ on $Y$ is $\eta$-cobounded, we have that $Y\subset (Hy)^{+\eta}$. Hence $Ko\subset (Hy)^{+\theta_2+\eta}$. In particular, for every $k\in K$, there exists $h_k\in H$ such that $|ko-h_ky|\le \theta_2+\eta$. Let $K'$ be a set of representatives of the set $H\backslash K$  of right cosets of $H$. Then the set $K''=\{h_k^{-1}k\colon k\in K'\}$ is a set of representatives of $H\backslash K$. We claim that $K''\subset U$. Let $k\in K'$. By the triangle inequality,
		$$|h_k^{-1}ko-o|=|ko-h_ko|\le |ko-h_ky|+|h_ky-h_ko|.$$
		Thus, $|h_k^{-1}ko-o|\le 2\theta_2+\eta+1$. This proves the claim. Consequently, 
		$$[K:H]\le |K''|\le |U|\le\theta.$$
	\end{enumerate}
\end{proof}

	\section{Constricting elements}
\label{sec:constricting-elements}

\begin{conv}
	In this section, we fix:
	\begin{itemize}[label=$\blacktriangleright$]
	\item Constants $\mu\ge 1$ and $\nu$, $\delta\ge 0$.
	\item A $(\mu,\nu)$-path system group $(G,X,\scrP)$.
	\item A $\delta$-constricting element $(g,A)$.
	\item A $\delta$-constricting map $\pi_A \colon X\to A$.
\end{itemize}
\end{conv}

\subsection{A \texorpdfstring{$G$}--invariant family}

The set of $G$-translates of $A$ is a $G$-invariant family of $\delta$-constricting subsets. Indeed, consider the stabilizer $\Stab(A)$ of $A$ and fix a set $R_g$ of representatives of $G/\Stab(A)$. Let $u\in G$ and $u_0\in R_g$ such that $uA=u_0A$. The map $\pi_{uA}\colon X\to uA$ defined as
$$\forall\,x\in X,\quad \pi_{uA}(x)=u_0\pi_A(u_0^{-1}x).$$ 
is then $\delta$-constricting since $\scrP$ is $G$-invariant. Moreover, the element $(ugu^{-1},uA)$ is $\delta$-constricting. To cope with the possible lack of $\langle ugu^{-1}\rangle$-equivariance of the map $\pi_{uA}\colon X\to uA$, we make the following observation:

\begin{prop}
	\label{prop:G-invariant}
	There exists $\theta\ge 0$ satisfying the following. Let $u\in G$. Then:
	\begin{enumerate}[label=(\roman*)]
		\item For every $x\in X$, we have $|\pi_{uA}(x)-u\pi_A(u^{-1}x)|\le \delta$.
		\item For every $Y\subset X$, we have $|\diam_{uA}(Y)-\diam(u\pi_A(u^{-1}Y))|\le \theta$.
	\end{enumerate}
\end{prop}

\begin{proof}
	Let $\theta_0=\theta_0(\delta)\ge 0$ be the constant of \autoref{prop:elementary properties}. We put $\theta=2\theta_0$. Let $u\in G$.
	\begin{enumerate}[label=(\roman*)]
		\item Let $x\in X$. Denote $y=u^{-1}x$. Let $u_0\in R_g$ such that $uA=u_0A$. We see that,
		$$|\pi_{uA}(x)-u\pi_A(u^{-1}x)|=|u_0\pi_{A}(u_0^{-1}x)-u\pi_A(u^{-1}x)|=|\pi_A(u_0^{-1}uy)-u_0^{-1}u\pi_A(y)|.$$
		Since $u_0^{-1}u\in Stab(A)$, it follows from \autoref{prop:elementary properties} (2) \emph{Coarse equivariance} that $|\pi_{uA}(x)-u\pi_A(u^{-1}x)|\le \theta_0$.
		\item Let $Y\subset X$. Let $y,y'\in Y$. By the triangle inequality,
		\begin{align*}
			\big||\pi_{uA}(y)-\pi_{uA}(y')|-&|u\pi_A(u^{-1}y)-u\pi_A(u^{-1}y')|\big|\le\\ &|\pi_{uA}(y)-u\pi_A(u^{-1}y)|+|u\pi_A(u^{-1}y')-\pi_{uA}(y')|.
		\end{align*}
		It follows from (i) that $$\max\qty{|u\pi_{uA}(y)-u\pi_A(u^{-1}y)|, |u\pi_A(u^{-1}y')-\pi_{uA}(y')|}\le \theta_0.$$
		Hence, we have $|\diam_{uA}(Y)-\diam(u\pi_A(u^{-1}Y))|\le 2\theta_0$.
	\end{enumerate}
	
\end{proof}

\subsection{Finding a constricting element}
The goal of this subsection is to combine \autoref{prop:alternative of intersection and hausdorff distance} and \autoref{prop:diameter}. We suggest to compare (ii) below with the property (BS2) of the buffering sequences.

\begin{prop}
	\label{prop:alternative of projection and hausdorff distance}
	Let $\eta\ge 0$. There exists $\theta \ge 1$ satisfying the following. Let $(H,Y)$ be an $\eta$-quasi-convex subgroup of $G$. Then:
	\begin{enumerate}[label=(\roman*)]
		\item For every $u\in G$, if $\diam_{uA}(Y)> \theta$, then $uA\subset Y^{+\theta}$.
		\item Let $H\le K\le G$. If $[K:H]> \theta$, then there exists $k\in K$ such that $\diam_{kA}(Y)\le \theta$.
	\end{enumerate}
\end{prop}

\begin{proof}
	
	Let $\eta\ge 0$. Let $\theta=\theta(\eta) \ge 1$. Its exact value will be precised below. It follows from \autoref{prop:elementary properties} (6) \emph{Morseness} and (7) \emph{Coarse invariance} that there exists $\theta_0 \ge 0$ such that the element $(g,A)$ is $\theta_0$-quasi-convex. Let $\theta_1=\max\{\eta,\theta_0\}$. By \autoref{prop:diameter}, there exist $\theta_2\ge 0$, $\zeta\ge 0$ depending both on $\theta_1$ such that for every $u\in G$ and for every $\theta_1$-quasi-convex subset $Y\subset X$, we have
	$$\diam_{uA}(Y)-\zeta \le \diam(uA^{+\theta_2}\cap Y)\le \diam_{uA}(Y)+\zeta.$$
	According to \autoref{prop:elementary properties} (6) \emph{Morseness} and (7) \emph{Coarse invariance}, there exist $\theta_3=\theta_3(\theta_2)\ge 0$ such that the element $(g,A^{+\theta_2})$ is $\theta_3$-quasi-convex. Let $\theta_4=\max\{\eta,\theta_3\}$. Let $\theta_5=\theta_5(\theta_4,g,A) \ge 1$ be the constant of \autoref{prop:alternative of intersection and hausdorff distance}. Finally, we put $\theta=\theta_5+\zeta$. Let $(H,Y)$ be an $\eta$-quasi-convex subgroup of $G$.
	
	\begin{enumerate}[label=(\roman*)]
		\item Let $u\in G$. Assume that $\diam_{uA}(Y)>\theta$. According to \autoref{prop:diameter}, we have $\diam(uA^{+\theta_2}\cap Y)>\theta_5$ and according to \autoref{prop:alternative of intersection and hausdorff distance} (i) this implies that $uA\subset Y^{+\theta_5}\subset Y^{+\theta}$.
		
		\item Let $H\le K \le G$. We argue by contraposition. Assume that for every $k\in K$, we have $\diam_{kA}(Y)> \theta$. According to \autoref{prop:diameter}, for every $k\in K$, we have $\diam(kA^{+\theta_2}\cap Y)> \theta_5$ and according to \autoref{prop:alternative of intersection and hausdorff distance} (ii) this implies that $[K:H]\le \theta_5\le \theta$.
	\end{enumerate}
	
\end{proof}

\subsection{Elementary closures}

The elementary closure of $(g,A)$ could be thought as the set of elements $u\in G$ such that $uA$ is ``parallel'' to $A$:
\begin{df}
	The \emph{elementary closure of} $(g,A)$ \emph{in} $G$ is defined as
	$$E(g,A)=\{u\in G\colon\, d_{\operatorname{\operatorname{Haus}}}(uA,A)<\infty\}.$$
	Observe that $E(g,A)$ is a subgroup of $G$ since $d_{\operatorname{Haus}}$ is a pseudo-distance. 
\end{df}
This subsection is devoted to provide a further description $E(g,A)$.  We suggest to compare the proposition below with the property (BS1) of the buffering sequences.
\begin{prop}
	\label{prop:zero-one constricting}
	There exists $\theta \ge 1$ satisfying the following:
	\begin{enumerate}[label=(\roman*)]
		\item For every $u\in G$, we have 
		$$\max\{\diam_{uA}(A),\diam_{A}(uA)\}>\theta \quad \Longleftrightarrow \quad d_{\operatorname{Haus}}(uA,A)\le\theta.$$
		\item $E(g,A)=\{u\in G\colon\, d_{\operatorname{Haus}}(uA,A)\le \theta\}$.
		\item $[E(g,A):\langle g \rangle]\le \theta$.
	\end{enumerate}
\end{prop}

\begin{proof}
	Let $\theta_0\ge 0$ be the constant of \autoref{prop:G-invariant}. According to \autoref{prop:elementary properties} (6) \emph{Morseness}, there exists $\theta_1\ge 0$ such that the element $(g,A)$ is $\theta_1$-quasi-convex. Let $\theta_2=\theta_2(\theta_1)\ge 1$ be the constant of \autoref{prop:alternative of projection and hausdorff distance}. We put $\theta=\theta_0+\theta_2$.
	
	\begin{cla}
		\label{lem:Hausdorff-Diameter}
		Let $u\in G$. If $d_{\operatorname{Haus}}(uA,A)<\infty$, then $\diam_{uA}(A)=\infty$. 
	\end{cla}
	
	Let $u\in G$. Assumme that $d_{\operatorname{Haus}}(uA,A)<\infty$ and denote $\varepsilon=d_{\operatorname{Haus}}(uA,A)+1$. By \autoref{prop:diameter}, there exist $\theta_3$, $\zeta\ge 0$ such that for every $u\in G$ we have
	$$\diam_{uA}(A)-\zeta \le \diam(uA^{+\theta_3}\cap A^{+\varepsilon})\le \diam_{uA}(A)+\zeta.$$
	Note that $uA\subset uA^{+\theta_3}\cap A^{+\varepsilon}$ and $\diam(uA)=\diam(A)$. Since the action of $G$ on $X$ is proper and since the element $g$ has infinite order, we have that $\diam(A)=\infty$. Consequently, we have $\diam(uA^{+\theta_3}\cap A^{+\varepsilon})=\infty$. Finally, it follows from \autoref{prop:diameter} that $\diam_{uA}(A)=\infty$. This proves the claim.
	
	\begin{enumerate}[label=(\roman*)]
		\item Let $u\in G$. Assume that $\max\{\diam_{uA}(A), \diam_{A}(uA)\}> \theta$. By \autoref{prop:G-invariant},
		$$\diam_{u^{-1}A}(A)\ge \diam_{A}(u^{-1}\pi_A(uA))-\theta_0.$$
		Hence, $\diam_{u^{-1}A}(A)>\theta_2$. It follows from \autoref{prop:alternative of projection and hausdorff distance} (i) that $uA\subset A^{+\theta}$ and  $u^{-1}A\subset A^{+\theta}$. Hence $d_{\operatorname{Haus}}(uA,A)\le \theta$. The converse follows from the claim above.
		
		\item This follows from (i) and the claim above.
		
		\item This follows from (i), (ii) and \autoref{prop:alternative of projection and hausdorff distance} (ii).
	\end{enumerate}
\end{proof}

Finally, we obtain an algebraic description of $E(g,A)$.

\begin{cor}
	\label{cor:elementary subgroup}
	There exist $\theta\ge 1$ and $M\in\zinterval{1}{\theta}$ such that for every $u\in G$, the following statements are equivalent:
	\begin{enumerate}[label=(\roman*)]
		\item $u\in E(g,A)$.
		\item There exists $p\in\{-1,1\}$ such that $ug^{M}u^{-1}=g^{p M}$.
		\item There exist $m,n\in \Z-\{0\}$ such that $ug^mu^{-1}=g^n$.
	\end{enumerate}
	Further, let $E^{+}(g,A)=\{u\in G\colon\, ug^{M}u^{-1}=g^{M}\}$. Then $[E(g,A):E^+(g,A)]\le 2$.
\end{cor}

\begin{proof}
	By \autoref{prop:zero-one constricting} $(ii)$, there exists $\theta_0\ge 1$ such that $[E(g,A):\langle g\rangle]\le \theta_0$. Let $\theta=\theta_0!$ We construct $M\in\zinterval{1}{\theta}$. First, we claim that there exists a subgroup $K\le \langle g\rangle$ such that $K\trianglelefteq E(g,A)$ and $[E(g,A):K]\le \theta$. Consider the natural action of $E(g,A)$ by right multiplication on the set $\langle g\rangle \backslash E(g,A)$ of right cosets of $\langle g\rangle$. This gives an homomorphism $\phi \colon E(g,A) \to \operatorname{Sym}(\langle g\rangle \backslash E(g,A))$. Choose $K=Ker(\phi)$. Note that $\langle g\rangle=\{h\in E(g,A)\colon \phi(h)(\langle g\rangle)\}=\langle g\rangle$. Thus, $K\le \langle g\rangle$. Morover, $K\trianglelefteq E(g,A)$. Further, we have that $|Sym(\langle g\rangle \backslash E(g,A))|=[E(g,A):\langle g\rangle]!$ and hence $[E(g,A):K]$ divides $[E(g,A):\langle g\rangle]!$ Therefore, $[E(g,A):K]\le \theta$. This proves the claim. Now, since the element $g$ has infinite order, the subgroup $E(g,A)$ is infinite. Hence, since $[E(g,A):K]<\infty$ there exists $M\ge 1$ such that  $K=\langle g^{M}\rangle$. Finally, we remark that $M$ is equal to the order of the element $\phi(g)$. Hence, $M\le \theta$.
	
	Let $u\in G$. The implication $(ii)\Rightarrow (iii)$ already holds.
	
	$(i)\Rightarrow (ii)$. Assume that $u\in E(g,A)$. Since the subgroup $\langle g^{M}\rangle$ is normal in $E(g,A)$, there exists $p \in\Z$ such that $ug^{M}u^{-1}=g^{p M}$. In particular, 
	$$\langle g^M\rangle=u\langle g^{M} \rangle u^{-1} = \langle ug^{M}u^{-1} \rangle=\langle g^{p M} \rangle.$$
	Hence, if $p\not\in\{-1,+1\}$, then $\langle g^M\rangle \not\subset \langle g^{p M}\rangle$. Contradiction.
	
	$(iii)\Rightarrow (i)$. Assume that there exist $m,n\in \Z-\{0\}$ such that $ug^mu^{-1}=g^n$. Since both $\langle g^m\rangle$ and $\langle g^n\rangle$ have finite index in $\langle g\rangle$, there exist $\zeta\ge 0$ the actions of $\langle ug^mu^{-1}\rangle $ on $uA$ and of $\langle g^n \rangle$ on $A$ are both $\zeta$-cobounded. Let $x\in uA$ and $y\in A$. We obtain $d_{\operatorname{Haus}}(uA,A)\le \zeta+|x-y|$. Hence $d_{\operatorname{Haus}}(uA,A)<\infty$.
	
	Finally, let $E^{+}(g,A)=\{u\in G\colon\, ug^{M}u^{-1}=g^{M}\}$. We prove that $[E(g,A):E^+(g,A)]\le 2$. It is enough to assume that $E(g,A)\neq E^+(g,A)$. Let $u,v\in E(g,A)-E^+(g,A)$. We show that $v^{-1}u\in E^+(g,A)$. Since $ug^{M}u^{-1}=vg^{M}v^{-1}=g^{-M}$, we have $v^{-1}ug^{M}u^{-1}v=v^{-1}g^{-M}v=g^{M}$ and therefore $v^{-1}u\in E^+(g,A)$. Hence $[E(g,A):E^+(g,A)]= 2$
\end{proof}

\subsection{Forcing a geometric separation}
In this subsection, we build large powers of our constricting element $(g,A)$ to produce a translate $Y'$ of a subset $Y$ so that the distance between their projections to a preferred $G$-translate of $A$ is large. We will do it in two different ways. We will apply these results to verify (BS4) in the construction of buffering sequences. Our main tool will be:

\begin{lem}
	\label{lem:separating projections}
	There exists $\theta\ge 0$ such that for every $x,x'\in X$ and for every $m\in \Z$,
	$$|x-g^mx'|_A\ge |m|\stlen{g}-|x-x'|_A-\theta.$$
\end{lem}

\begin{proof}
	Let $\theta=\theta(\delta)\ge 0$ be the constant of \autoref{prop:elementary properties}. Let $x, x' \in X$. Let $m\in \Z$. If $m=0$, then there is nothing to do. Assume that $m\neq 0$. By the triangle inequality,
	$$|x-g^mx'|_A\ge |\pi_A(x)-g^m\pi_A(x)|-|x-x'|_A-|g^m\pi_A(x')-\pi_A(g^mx')|.$$
	Note that
	$$\frac{1}{|m|}|\pi_A(x)-g^m\pi_A(x)|\ge \inf_{n\ge 1} \frac{1}{n}|\pi_A(x)-g^n\pi_A(x)|= \stlen{g}.$$
	By \autoref{prop:elementary properties} (2) \emph{Coarse equivariance}, we have $|g^m\pi_A(x')-\pi_A(g^mx')|\le \theta$. Therefore, we have
	$|x-g^mx'|_A\ge |m|\stlen{g}-|x-x'|_A-\theta$.
\end{proof}

The first way of forcing a geometric separation will be applied to the study of the relative exponential growth rates:
\begin{prop}
	\label{prop:geometric separation 1}
	For every $\varepsilon$, $\theta\ge 0$, there exists $M\ge 1$ with the following property. Let $H\le G$ be a subgroup. Let $Y\subset X$ be an $H$-invariant subset. If $\diam_A(Y)\le \varepsilon$, then for every $u\in \langle g^{M}, H\cap E(g,A)\rangle-H\cap E(g,A)$, we have  $d_{A}(Y,uY)>\theta$.
\end{prop}

\begin{proof}
	Let $\varepsilon$, $\theta\ge 0$. Let $\theta_0\ge 0$ be the constant of \autoref{prop:elementary properties}. By \autoref{lem:separating projections}, there exists $\theta_1\ge 0$ such that for every $x,x'\in X$ and for every $m\in \Z$,
	$$|x-g^mx'|_A\ge |m|\stlen{g}-|x-x'|_A-\theta_1.$$
	Combining \autoref{lem:qi embedding of quasiconvex subgroups} and \autoref{prop:elementary properties} (6) \emph{Morseness}, we obtain $\stlen{g}>0$. According to \autoref{cor:elementary subgroup}, there exists $M_0\ge 1$ such that 
	$$E(g,A)=\left\{u\in G\colon\,\exists\,p\in\{-1,+1\}\, ug^{M_0}u^{-1}=g^{pM_0}\right\}.$$ Let $m_0>\frac{\theta-2\varepsilon-2\theta_0-\theta_1}{M_0\stlen{g}}$. We put $M=M_0m_0$.  
	
	Let $H\le G$ be a subgroup. Let $Y\subset X$ be an $H$-invariant subset. Assume that $\diam_A(Y)\le \varepsilon$. Let $u\in \langle g^{M}, H\cap E(g,A)\rangle-H\cap E(g,A)$ and $y,y'\in Y$. It follows from \autoref{cor:elementary subgroup} that there exists $n\in \Z$ multiple of $M$ and $f\in H\cap E(g,A)$ such that $u=g^{n}f$. By the triangle inequality,
	
	\begin{align*}
		|y-g^{n}fy'|_A\ge |y-g^{n}y'|_A-|\pi_A(g^{n}y')-g^{n}\pi_A(y')|
		-|y'-fy'|_A-|g^{n}\pi_A(fy')-\pi_A(g^{n}fy')|.
	\end{align*}
	By \autoref{lem:separating projections},
	$$|y-g^{n}y'|_A\ge |n|\stlen{g}-|y-y'|_A-\theta_1$$
	Note that $u\not\in H\cap E(g,A)$ implies $n\neq 0$. Hence $|n|\ge |M|$. Since $f\in H$ and $\diam_A(Y)\le \varepsilon$,
	$$\max\{|y-y'|_A, |y'-fy'|_A\}\le \varepsilon.$$
	By \autoref{prop:elementary properties} (2) \emph{Coarse equivariance}, 
	$$\max\{|\pi_A(g^{n}y')-g^{n}\pi_A(y')|,|g^{n}\pi_A(fy')-\pi_A(g^{n}fy')|\}\le \theta_0.$$
	Since the elements $y,y'$ are arbitrary, we obtain $d_{A}(Y,uY)>\theta$.
\end{proof}
The second way of forcing a geometric separation will be applied to the study of the quotient exponential growth rates:
\begin{prop}
	\label{prop:geometric separation 2}
	For every $\varepsilon$, $\theta\ge 0$, there exist $M\ge 1$ and $f\colon G\times X \to \{1_G,g^M\}$ with the following property. Let $Y\subset X$ be subset. If $\diam_A(Y)\le \varepsilon$, then for every $u\in G$ and for every $y\in Y$, we have $d_{uA}(y,uf(u,y)Y)>\theta$.
\end{prop}

\begin{proof}
	Let $\varepsilon$, $\theta\ge 0$. Let $\theta_0\ge 0$ be the constant of \autoref{prop:G-invariant}. By \autoref{lem:separating projections}, there exists $\theta_1\ge 0$ such that for every $x,x'\in X$ and for every $m\in \Z$,
	$$|x-g^mx'|_A\ge |m|\stlen{g}-|x-x'|_A-\theta_1.$$
	Combining \autoref{lem:qi embedding of quasiconvex subgroups} and \autoref{prop:elementary properties} (6) \emph{Morseness}, we obtain $\stlen{g}>0$. We put $M>\frac{2\theta+2\varepsilon+8\theta_0+\theta_1}{\stlen{g}}$. Then, for every $u\in G$ and for every $x\in X$, there exists $f(u,x)\in\{1_G,g^M\}$ such that $|u^{-1}x-f(u,x)|_A>\theta+\varepsilon+4\theta_0$: if $|u^{-1}x-x|_A>\theta+\varepsilon+4\theta_0$, we choose $f(u,x)=1_G$, otherwise we choose $f(u,x)=g^M$.
	This defines $f\colon G\times X\to \{1_G,g^M\}$. 
	
	Let $Y\subset X$ be a subset. Assume that $\diam_A(Y)\le \varepsilon$. Let $u\in G$. Let $y,y'\in Y$. By abuse of notation, we write $f$ instead of $f(u,y)$. By the triangle inequality,
	\begin{align*}
		|y-ufy'|_{uA}&\ge |y-ufy|_{uA}-|ufy-ufy'|_{uA},\\
		|y-ufy|_{uA}&\ge |u^{-1}y-fy|_A-|\pi_{uA}(y)-u\pi_A(u^{-1}y)|-|\pi_{uA}(ufy)-u\pi_A(fy)|,\\ 
		|ufy-ufy'|_{uA} &\le |\pi_{uA}(ufy)-uf\pi_A(y)|+|y-y'|_A+|uf\pi_A(y')- \pi_{uA}(ufy')|.
	\end{align*}
	By hypothesis, $|u^{-1}y-fy|_A>\theta+\varepsilon+4\theta_0$ and $|y-y'|_A\le \diam_A(Y)\le \varepsilon$. By \autoref{prop:G-invariant},
	$$\max\{|\pi_{uA}(y)-u\pi_A(u^{-1}y)|,|\pi_{uA}(ufy)-u\pi_A(fy)|\}\le \theta_0.$$
	$$\max\{|\pi_{uA}(ufy)-uf\pi_A(y)|, |uf\pi_A(y')- \pi_{uA}(ufy')|\}\le \theta_0.$$
	Since the element $y'$ is arbitrary, we obtain $d_{uA}(y,ufY)>\theta$.
\end{proof}

	\section{Growth of quasi-convex subgroups}
\label{sec:growth}

The goal of this section is to prove \autoref{th:entropy1} and \autoref{th:entropy2}.

\begin{conv}
	In this section, we fix:
	\begin{itemize}[label=$\blacktriangleright$]
	\item Constants $\mu\ge 1$ and $\nu$, $\delta$, $\eta\ge 0$.
	\item A $(\mu,\nu)$-path system group $(G,X,\scrP)$.
	\item A $\delta$-constricting element $(g_0,A_0)$.
	\item An infinite index $\eta$-quasi-convex subgroup $(H,Y)$ of $G$.
\end{itemize}
\end{conv}

We are going to replace the axis $A_0$ for $A_0'=E(g_0,A_0) A_0$. Note that $d_{\operatorname{Haus}}(A_0,A_0')<\infty$ (\autoref{prop:zero-one constricting} (ii)). Up to replacing $\delta$ for a larger constant, it follows from \autoref{prop:elementary properties} (7) \emph{Coarse invariance } and \autoref{cor:elementary subgroup} that the element $(g_0,A_0')$ is $\delta$-constricting. By abuse of notation, we still denote $A_0=A_0'$. In this new setting, we have $kA_0=A_0$, for every $k\in E(g_0,A_0)$. Let $\theta_0=\theta_0(\delta,\eta)\ge 1$ be the constant of \autoref{prop:alternative of projection and hausdorff distance}. Since $[G:H]=\infty$, there exist $u\in G$ such that $\diam_{uA_0}(Y)\le \theta_0$  (\autoref{prop:alternative of projection and hausdorff distance} (ii)). We denote $(g,A)=(ug_0u^{-1},uA_0)$.

\subsection{Case $\omega(H)<\omega(G)$}

In this subsection we prove:

\begin{thm}[\autoref{th:entropy1}]
\label{th:entropy1-body}
Assume that
	\begin{enumerate}[label=(\roman*)]
	\item $\omega(H)<\infty$.
	\item The action of $H$ on $X$ is divergent.
\end{enumerate}	
Then $\omega(H)<\omega(G)$.
\end{thm}

We require the following.

\begin{prop}[\autoref{th:amalgam}]
	\label{prop:injectivity}
	There exist $M\ge 1$ satisfying the following:
	\begin{enumerate}[label=(\roman*)]
		\item $E(g,A)$ is a finite extension of $\group{g}$.
		\item $H\cap E(g,A)$ is a finite proper subgroup of $\langle g^{M}, H\cap E(g,A)\rangle$. 
		\item The natural homomorphism $H\ast_{H\cap E(g,A)} \langle g^{M}, H\cap E(g,A)\rangle\to G$ is injective.
	\end{enumerate}
\end{prop}

\begin{proof}
 	The subgroup $E(g,A)$ is a finite extension of $\langle g \rangle$ (\autoref{prop:zero-one constricting} (iii)). This proves (i). Since $\diam_A(Y)\le \theta_0$ and the action of $H\cap E(g,A)$ on $Y \cap A^{+\rho}$ for $\rho=d(A,Y)$ is proper and cobounded, the subgroup $H\cap E(g,A)$ is finite (\autoref{prop:diameter}). Further, since $g$ has infinite order, $H\cap E(g,A)$ must be a proper subgroup of $\langle g^{M}, H\cap E(g,A)\rangle$. This proves (ii). 
	
	The rest of the proof is devoted to establish (iii). Let $\theta_1=\theta_1(\delta)\ge 0$ be the constant of \autoref{prop:G-invariant}. Let $\varepsilon=\max\{\theta_0+2\theta_1,d(A,Y)\}$. Let $L=L(\delta,\varepsilon,0)\ge 0$ be the constant of \autoref{cor:buffer sequence}. By \autoref{prop:geometric separation 1}, there exists $M\ge 1$ such that for every $u\in \langle  g^{M} ,H\cap E(g,A)\rangle- H\cap E(g,A)$, we have $d_{A}(Y,uY)>L-2\theta_1$. Let $\phi\colon H\ast_{H\cap E(g,A)} \langle g^{M} ,H\cap E(g,A)\rangle\to G$ be the natural homomorphism. Let $w\in H\ast_{H\cap E(g,A)} \langle g^{M} ,H\cap E(g,A)\rangle$ such that $w\neq 1$. We are going to prove that $\phi(w)\neq 1$. Note that the homomorphisms $\phi_{|H}$ and $\phi_{| \langle g^{M} ,H\cap E(g,A)\rangle}$ are injective. If $w\in H\cup \langle g^{M} ,H\cap E(g,A)\rangle$, then $\phi(w)\neq 1$. Assume that $w\not\in H\cup \langle g^{M} ,H\cap E(g,A)\rangle$. Note that if there exists a conjugate $w'$ of $w$ such that $\phi(w')\neq 1$, then $\phi(w)\neq 1$. Up to replacing $w$ by a cyclic conjugate, there exist $n\ge 1$ and a sequence $h_1,k_1,\cdots, h_n,k_n\in G$ such that $w=h_1k_1\cdots h_nk_n$ and such that for every $i\in\{1,\cdots,n\}$ we have $h_i\in  H-H\cap E(g,A)$ and $k_i\in \langle g^{M} ,H\cap E(g,A)\rangle-H\cap E(g,A)$. For every $i\in\zinterval{1}{n}$, we denote $u_i=h_1k_1\cdots h_{i}$ and $v_i=h_1k_1\cdots h_{i}k_{i}$. We also denote $v_0=1_G$. 
	
	We are going to prove that the sequence $v_0Y,u_1A,v_1Y,\cdots, u_nA,v_nY$ is $(\delta,\varepsilon,L)$-buffering on $\{u_iA\}$ and then apply \autoref{cor:buffer sequence}. Let $i\in\zinterval{1}{n}$. Let us prove (BS1). Assume for a moment that $i\neq n$.  Since we had modified the axis $A_0$ above, for every $j\in\zinterval{1}{n}$, we have $k_jA=A$. Hence
	\begin{align*}
		\pi_{u_iA}(u_{i+1}A)&=\pi_{v_iA}(u_{i+1}A),\\
		\pi_{u_{i+1}A}(u_iA)&=\pi_{u_{i+1}A}(v_iA).
	\end{align*}
	By \autoref{prop:G-invariant},
	\begin{align*}
		\diam_{v_iA}(u_{i+1}A)&\le \diam(v_i\pi_A(h_iA))+\theta_1,\\
		\diam_{u_{i+1}A}(v_iA)&\le \diam(u_{i+1}\pi_A(h_i^{-1}A))+\theta_1,\\
		\diam_A(h_i^{-1}A)&\le \diam_{h_iA}(A)+\theta_1.
	\end{align*}
	By \autoref{prop:zero-one constricting} (i) and (ii), for every $u\not \in E(g,A)$, we have $\max\{\diam_A(uA),\diam_{uA}(A)\}\le \theta_0$. Consequently, $$\max\{\diam_{u_iA}(u_{i+1}A),\diam_{u_{i+1}A}(u_iA)\}\le \theta_0+2\theta_1 \le \varepsilon.$$
	Let us prove (BS2). Note that,
	\begin{align*}
		\pi_{u_iA}(v_{i-1}Y)&=\pi_{u_iA}(u_iY),\\
		\pi_{u_{i}A}(v_iY)&=\pi_{v_{i}A}(v_{i}Y).
	\end{align*}
	By \autoref{prop:G-invariant},
	\begin{align*}
		\diam_{u_iA}(u_iY)&\le \diam(u_i\pi_A(Y))+\theta_1,\\
		\diam_{v_{i}A}(v_{i}Y)&\le \diam(v_{i}\pi_A(Y))+\theta_1.
	\end{align*}
	Since $\diam_A(Y)\le \theta_0$, we obtain $$\max\{\diam_{u_iA}(v_{i-1}Y),\diam_{u_{i}A}(v_iY)\}\le \theta_0+\theta_1 \le \varepsilon.$$
	Let us prove (BS3). We have,
	$$\max\{d(u_iA,v_{i-1}Y),d(u_iA,v_{i}Y)\}=\max\{d(u_iA,u_iY),d(v_iA,v_iY)\}\le d(A,Y)\le \varepsilon.$$
	Let us prove (BS4). It follows from \autoref{prop:G-invariant} (i) that,
	$$d_{u_iA}(v_{i-1}Y,v_iY)\ge d_A(Y,k_iY)-2\theta_1.$$
	By the choice of $M$, we have $d_A(Y,k_iY)>L+2\theta_1$. Hence, we have $d_{u_iA}(v_{i-1}Y,v_iY)\ge L$. This proves that the sequence $v_0Y,u_1A,v_1Y,\cdots, u_nA,v_nY$ is $(\delta,\varepsilon,L)$-buffering on $\{u_iA\}$. It follows from \autoref{cor:buffer sequence} that $d_{u_{n}A}(Y,\phi(w)Y)>0$. Hence, $\phi(w)\neq 1$.
	
\end{proof}

\begin{proof}[Proof of \autoref{th:entropy1-body}]
\autoref{th:entropy1-body} is an immediate consequence of \autoref{prop:Dalbo} and \autoref{prop:injectivity}.
\end{proof}

\subsection{Case $\omega(G/H)=\omega(G)$}
In this subsection we prove:
\begin{thm}[\autoref{th:entropy2}]
	\label{th:entropy2-body}
	$\omega(G/H)=\omega(G)$.
\end{thm}

Recall that given $\phi\colon G\to G$, we say that  $G$ is $\phi$-\emph{coarsely} $G/H$ if there exist $\theta\ge 0$, $x\in X$ satisfying the following conditions:
\begin{enumerate}[label=(CQ\arabic*)]
	\item For every $u,v\in G$, if $\phi(u)H=\phi(v)H$, then $|\phi(u)x-\phi(v)x|\le \theta$.
	\item For every $u\in G$, $|ux-\phi(u)x|\le \theta$.
\end{enumerate}

We require the following.
\begin{prop}
	\label{prop:geometric separation of projections}
	There exist $M\ge 1$ and a map $f\colon G \to \{1_G,g^M\}$ with the following property. Let $\phi\colon G\to G$, $u\mapsto uf_u$. Then $G$ is $\phi$-coarsely $G/H$.
\end{prop}

We prove some preliminar lemmas.

\begin{lem}
	\label{lem:quotient 1}
	There exists $\theta\ge 0$ such that for every $m\in \Z$, we have $\diam_A(g^{m}Y)\le \theta$.
\end{lem}

\begin{proof}
	Let $\theta_1\ge 0$ be the constant of \autoref{prop:elementary properties}. We put $\theta=\theta_0+2\theta_1$. Let $m\in \Z$. Let $x,x'\in Y$. By the triangle inequality,
	$$|g^{m}x-g^{m}x'|_A\le |\pi_A(g^{m}x)-g^{m}\pi_A(x)|+|x-x'|_A+|g^{m}\pi_A(x')-\pi_A(g^{m}x')|.$$
	By \autoref{prop:elementary properties} (2) \emph{Coarse equivariance},
	$$\max\{|\pi_A(g^{m}x)-g^{m}\pi_A(x)|,|g^{m}\pi_A(x')-\pi_A(g^{m}x')|\}\le \theta_1.$$
	Moreover, we have $|x-x'|_A\le \diam_A(Y)\le \theta_0$. Since $x,x'$ are arbitrary, we obtain $\diam_A(g^{m}Y)\le \theta_0+2\theta_1$.
\end{proof}

\begin{lem}
	\label{lem:quotient 2}
	For every $\varepsilon\ge 0$, there exists $\theta \ge 0$ with the following property. Let $A_1,A_2\subset X$ be $\delta$-constricting subsets such that $d_{\operatorname{Haus}}(A_1,A_2)\le \varepsilon$. Let $x\in A_1^{+\varepsilon}$ and $y\in A_2^{+\varepsilon}$ such that $|x-y|_{A_1}\le \varepsilon$. Then $|x-y|\le \theta$.
\end{lem}

\begin{proof}
	Let $\theta_1\ge 0$ be the constant of \autoref{prop:elementary properties}. Let $\varepsilon\ge 0$. Let $\theta\ge 0$. Its exact value will be precised below. Let $A_1,A_2\subset X$ be $\delta$-constricting subsets such that $d_{\operatorname{Haus}}(A_1,A_2)\le \varepsilon$. Let $x\in A_1^{+\varepsilon}$ and $y\in A_2^{+\varepsilon}$ such that $|x-y|_{A_1}\le \varepsilon$. By the triangle inequality,
	$$|x-y|\le  |x-\pi_{A_1}(x)|+|x-y|_{A_1}+|\pi_{A_1}(y)-y|.$$
	Since $x,y\in A_1^{+2\varepsilon+1}$, it follows from \autoref{prop:elementary properties} (1) \emph{Coarse nearest-point projection} that
	$$\max\{|x-\pi_{A_1}(x)|,|\pi_{A_1}(y)-y|\}\le \mu(2\varepsilon+1)+\theta_1.$$
	Finally, we put $\theta=\varepsilon+2\mu(2\varepsilon+1)+2\theta_1$.
\end{proof}

We are ready to prove \autoref{prop:geometric separation of projections}:

\begin{proof}[Proof of \autoref{prop:geometric separation of projections}]
	Let $\theta_1\ge 0$ be the constant of \autoref{prop:G-invariant}. Let $\theta_2\ge0$ be the constant of \autoref{prop:zero-one constricting}. Let $\theta_3\ge 0$ be the constant of \autoref{lem:quotient 1}. Let $\varepsilon=\max\{\theta_2+2\theta_1,\theta_1+\theta_3, d(A,Y)+1\}$. In particular, there exists $y\in A^{+\varepsilon}\cap Y$. Let $\theta_4=\theta_4(\delta,\varepsilon)\ge 0$ be the constant of \autoref{prop:Behrstock's inequality 2}. By \autoref{prop:geometric separation 2}, there exist $M\ge 1$ and $f\colon G\to \{1_G,g^M\}$ such that for every $u\in G$, we have $d_{uA}(y,uf(u)Y)>\theta_4$. For every $u\in G$, we denote $f_u=f(u)$ and we put $\phi\colon G\to G, u\mapsto uf_u$. Let $\theta_5=\theta_5(\varepsilon)\ge 0$ be the constant of \autoref{lem:quotient 2}. We put $\theta=\max\{|y-g^My|,\theta_5\}$. We are going to prove that $G$ is $\phi$-coarsely $G/H$ with respect to $y$ and $\theta$.
	
	In order to prove (CQ1), we just need to observe that for every $u\in G$, we have 
	$$|uy-uf_uy|=|y-f_uy|\le |y-g^My|\le \theta.$$
	Let us prove (CQ2). Let $u,v\in G$. Assume that $uf_uH=vf_vH$. We claim that $d_{\operatorname{Haus}}(uA,vA)\le \theta_2$. By \autoref{prop:zero-one constricting} (i), it suffices to prove that
	$$\max\{\diam_{v^{-1}uA}(A),\diam_A(v^{-1}uA)\}>\theta_2.$$
	We argue by contradiction. Assume instead that $\max\{\diam_{v^{-1}uA}(A),\diam_A(v^{-1}uA)\}\le \theta_2$. We are going to prove that the sequence $uA, uf_uY, vA$ is $(\delta,\varepsilon,0)$-buffering on $\{uA,vA\}$ and then apply \autoref{prop:Behrstock's inequality 2}. Note that the condition (BS4) is void in this case. Let us prove (BS1). By \autoref{prop:G-invariant},
	\begin{equation*}
		\begin{aligned}
			\diam_{uA}(vA)&\le \diam(u\pi_{A}(u^{-1}vA))+\theta_1,\\ \diam_{vA}(uA)&\le \diam(v\pi_A(v^{-1}uA))+\theta_1,\\
			\diam_A(u^{-1}vA)&\le \diam_{v^{-1}uA}(A)+\theta_1.
		\end{aligned}
	\end{equation*}
	Hence,
	$$\max\{\diam_{uA}(vA),\diam_{vA}(uA)\}\le \theta_2+2\theta_1\le \varepsilon.$$ 
	Let us prove (BS2). By \autoref{prop:G-invariant},
	\begin{equation*}
		\begin{aligned}
			\diam_{uA}(uf_uY)&\le \diam(u\pi_A(f_uY))+\theta_1,\\ \diam_{vA}(vf_vY)&\le \diam(v\pi_A(f_vY))+\theta_1.\\
		\end{aligned}
	\end{equation*}
	By \autoref{lem:quotient 1}, we have $\max\{\diam_A(f_uY),\diam_A(f_vY)\}\le \theta_3$. Hence, 
	$$\max\{\diam_{uA}(uf_uY),\diam_{vA}(vf_vY)\}\le \theta_1+\theta_3\le \varepsilon.$$ 
	Let us prove (BS3). The hypothesis $uf_uH=vf_vH$ implies $uf_uY=vf_vY$ and therefore
	$$\max\{d(uA,uf_uY),d(vA,uf_uY)\}=\max\{d(uA,uf_uY),d(vA,vf_vY)\}=d(A,Y)\le \varepsilon.$$ 
	Hence, the sequence $uA,uf_uY, vA$ is $(\delta,\varepsilon,0)$-buffering on $\{uA,vA\}$.  It follows from \autoref{prop:Behrstock's inequality 2} that
	$$\min\qty{d_{uA}(y,uf_uY),d_{vA}(y,uf_uY)}\le \theta_4.$$ 
	However, by construction, $$\min\qty{d_{uA}(y,uf_uY),d_{vA}(y,uf_uY)}> \theta_4.$$ 
	Contradiction. Therefore, $d_{\operatorname{Haus}}(uA,vA)\le \theta_2$. This proves the claim. In particular, $d_{\operatorname{Haus}}(uA,vA)\le \varepsilon$. Since $y\in A^{+\varepsilon}$, we have $uf_uy\in uA^{+\varepsilon}$ and $vf_vy\in vA^{+\varepsilon}$. Since $uf_uy,vf_vy\in uf_uY$, we have $|uf_uy-vf_vy|_{uA}\le \diam_{uA}(uf_uY) \le \varepsilon$. According to \autoref{lem:quotient 2}, $|uf_uy-vf_vy|\le \theta$. This proves (CQ2).
\end{proof}

\begin{proof}[Proof of \autoref{th:entropy2-body}]
\autoref{th:entropy2-body} is an immediate consequence of \autoref{prop:quasi-bijection} and \autoref{prop:geometric separation of projections}.
\end{proof}

	% References
	%\bibliographystyle{alpha}
	\bibliographystyle{abbrv} 
	\bibliography{biblio}

@article{abbott_property_2019,
	title = {Property  {P}-naive for acylindrically hyperbolic groups},
	volume = {291},
	issn = {0025-5874},
	url = {https://doi.org/10.1007/s00209-018-2094-1},
	doi = {10.1007/s00209-018-2094-1},
	number = {1-2},
	journal = {Mathematische Zeitschrift},
	author = {Abbott, Carolyn R. and Dahmani, François},
	year = {2019},
	mrnumber = {3936081},
	pages = {555--568},
}

@article{antolin_counting_2021,
	title = {Counting subgraphs in fftp graphs with symmetry},
	volume = {170},
	issn = {0305-0041},
	url = {https://doi.org/10.1017/S0305004119000422},
	doi = {10.1017/S0305004119000422},
	number = {2},
	journal = {Mathematical Proceedings of the Cambridge Philosophical Society},
	author = {Antolín, Yago},
	year = {2021},
	mrnumber = {4222436},
	pages = {327--353},
}

@article{arzhantseva_quasiconvex_2001,
	title = {On quasiconvex subgroups of word hyperbolic groups},
	volume = {87},
	issn = {0046-5755},
	url = {https://doi.org/10.1023/A:1012040207144},
	doi = {10.1023/A:1012040207144},
	number = {1-3},
	journal = {Geometriae Dedicata},
	author = {Arzhantseva, Goulnara N.},
	year = {2001},
	mrnumber = {1866849},
	pages = {191--208},
}

@article{arzhantseva_growth_2015,
	title = {Growth tight actions},
	volume = {278},
	issn = {0030-8730},
	url = {https://doi.org/10.2140/pjm.2015.278.1},
	doi = {10.2140/pjm.2015.278.1},
	number = {1},
	journal = {Pacific Journal of Mathematics},
	author = {Arzhantseva, Goulnara N. and Cashen, Christopher H. and Tao, Jing},
	year = {2015},
	mrnumber = {3404665},
	pages = {1--49},
}

@article{bestvina_characterization_2009,
	title = {A characterization of higher rank symmetric spaces via bounded cohomology},
	volume = {19},
	issn = {1016-443X},
	url = {https://doi.org/10.1007/s00039-009-0717-8},
	doi = {10.1007/s00039-009-0717-8},
	number = {1},
	journal = {Geometric and Functional Analysis},
	author = {Bestvina, Mladen and Fujiwara, Koji},
	year = {2009},
	mrnumber = {2507218},
	pages = {11--40},
}

@article{behrstock_hierarchically_2017,
	title = {Hierarchically hyperbolic spaces {I}: {Curve} complexes for cubical groups},
	volume = {21},
	issn = {1465-3060},
	url = {https://doi.org/10.2140/gt.2017.21.1731},
	doi = {10.2140/gt.2017.21.1731},
	number = {3},
	journal = {Geometry \& Topology},
	author = {Behrstock, Jason and Hagen, Mark F. and Sisto, Alessandro},
	year = {2017},
	mrnumber = {3650081},
	pages = {1731--1804},
}

@article{behrstock_hierarchically_2019,
	title = {Hierarchically hyperbolic spaces {II}: {Combination} theorems and the distance formula},
	volume = {299},
	issn = {0030-8730},
	url = {https://doi.org/10.2140/pjm.2019.299.257},
	doi = {10.2140/pjm.2019.299.257},
	number = {2},
	journal = {Pacific Journal of Mathematics},
	author = {Behrstock, Jason and Hagen, Mark F. and Sisto, Alessandro},
	year = {2019},
	mrnumber = {3956144},
	pages = {257--338},
}

@misc{behrstock_quasiflats_2020,
	title = {Quasiflats in hierarchically hyperbolic spaces},
	author = {Behrstock, Jason and Hagen, Mark F. and Sisto, Alessandro},
	year = {2020},
	note = {\_eprint: 1704.04271},
}

@article{behrstock_asymptotic_2006,
	title = {Asymptotic geometry of the mapping class group and {Teichmüller} space},
	volume = {10},
	issn = {1465-3060},
	url = {https://doi.org/10.2140/gt.2006.10.1523},
	doi = {10.2140/gt.2006.10.1523},
	journal = {Geometry and Topology},
	author = {Behrstock, Jason},
	year = {2006},
	mrnumber = {2255505},
	pages = {1523--1578},
}

@article{cashen_morse_2020,
	title = {Morse subsets of {CAT}(0) spaces are strongly contracting},
	volume = {204},
	issn = {0046-5755},
	url = {https://doi.org/10.1007/s10711-019-00457-x},
	doi = {10.1007/s10711-019-00457-x},
	journal = {Geometriae Dedicata},
	author = {Cashen, Christopher H.},
	year = {2020},
	mrnumber = {4056705},
	pages = {311--314},
}

@article{coulon_growth_2018,
	title = {Growth gap in hyperbolic groups and amenability},
	volume = {28},
	issn = {1016-443X},
	url = {https://doi.org/10.1007/s00039-018-0459-6},
	doi = {10.1007/s00039-018-0459-6},
	number = {5},
	journal = {Geometric and Functional Analysis},
	author = {Coulon, Rémi and Dal'Bo, Françoise and Sambusetti, Andrea},
	year = {2018},
	mrnumber = {3856793},
	pages = {1260--1320},
}

@book{coornaert_geometrie_1990,
	series = {Lecture {Notes} in {Mathematics}},
	title = {Géométrie et théorie des groupes},
	volume = {1441},
	isbn = {3-540-52977-2},
	publisher = {Springer-Verlag, Berlin},
	author = {Coornaert, Michel and Delzant, Thomas and Papadopoulos, Athanase},
	year = {1990},
	mrnumber = {1075994},
}

@article{dalbo_growth_2011,
	title = {On the growth of quotients of {Kleinian} groups},
	volume = {31},
	issn = {0143-3857},
	url = {https://doi.org/10.1017/S0143385710000131},
	doi = {10.1017/S0143385710000131},
	number = {3},
	journal = {Ergodic Theory and Dynamical Systems},
	author = {Dal'Bo, Françoise and Peigné, Marc and Picaud, Jean-Claude and Sambusetti, Andrea},
	year = {2011},
	mrnumber = {2794950},
	pages = {835--851},
}

@article{dahmani_growth_2019,
	title = {Growth of quasiconvex subgroups},
	volume = {167},
	issn = {0305-0041},
	url = {https://doi.org/10.1017/s0305004118000440},
	doi = {10.1017/s0305004118000440},
	number = {3},
	journal = {Mathematical Proceedings of the Cambridge Philosophical Society},
	author = {Dahmani, François and Futer, David and Wise, Daniel T.},
	year = {2019},
	mrnumber = {4015648},
	pages = {505--530},
}

@article{durham_boundaries_2017,
	title = {Boundaries and automorphisms of hierarchically hyperbolic spaces},
	volume = {21},
	issn = {1465-3060},
	url = {https://doi.org/10.2140/gt.2017.21.3659},
	doi = {10.2140/gt.2017.21.3659},
	number = {6},
	journal = {Geometry \& Topology},
	author = {Durham, Matthew Gentry and Hagen, Mark F. and Sisto, Alessandro},
	year = {2017},
	mrnumber = {3693574},
	pages = {3659--3758},
}

@article{gitik_growth_2020,
	title = {On growth of double cosets in hyperbolic groups},
	volume = {30},
	issn = {0218-1967},
	url = {https://doi.org/10.1142/S0218196720500356},
	doi = {10.1142/S0218196720500356},
	number = {6},
	journal = {International Journal of Algebra and Computation},
	author = {Gitik, Rita and Rips, Eliyahu},
	year = {2020},
	mrnumber = {4155416},
	pages = {1161--1166},
}

@incollection{gromov_hyperbolic_1987,
	series = {Math. {Sci}. {Res}. {Inst}. {Publ}.},
	title = {Hyperbolic groups},
	volume = {8},
	url = {https://doi.org/10.1007/978-1-4613-9586-7_3},
	booktitle = {Essays in group theory},
	publisher = {Springer, New York},
	author = {Gromov, M.},
	year = {1987},
	mrnumber = {919829},
	doi = {10.1007/978-1-4613-9586-7_3},
	pages = {75--263},
}

@article{gruber_infinitely_2018,
	title = {Infinitely presented graphical small cancellation groups are acylindrically hyperbolic},
	volume = {68},
	issn = {0373-0956},
	url = {http://aif.cedram.org/item?id=AIF_2018__68_6_2501_0},
	number = {6},
	journal = {Université de Grenoble. Annales de l'Institut Fourier},
	author = {Gruber, Dominik and Sisto, Alessandro},
	year = {2018},
	mrnumber = {3897973},
	pages = {2501--2552},
}

@book{de_la_harpe_topics_2000,
	series = {Chicago {Lectures} in {Mathematics}},
	title = {Topics in geometric group theory},
	isbn = {0-226-31719-6 0-226-31721-8},
	publisher = {University of Chicago Press, Chicago, IL},
	author = {de la Harpe, Pierre},
	year = {2000},
	mrnumber = {1786869},
}

@article{kapovich_nonamenability_2002,
	title = {The nonamenability of {Schreier} graphs for infinite index quasiconvex subgroups of hyperbolic groups},
	volume = {48},
	issn = {0013-8584},
	number = {3-4},
	journal = {L'Enseignement Mathématique. Revue Internationale. 2e Série},
	author = {Kapovich, Ilya},
	year = {2002},
	mrnumber = {1955607},
	pages = {359--375},
}

@misc{kapovich_geometry_2002,
	title = {The geometry of relative {Cayley} graphs for subgroups of hyperbolic groups},
	copyright = {Assumed arXiv.org perpetual, non-exclusive license to distribute this article for submissions made before January 2004},
	url = {https://arxiv.org/abs/math/0201045},
	publisher = {arXiv},
	author = {Kapovich, Ilya},
	year = {2002},
	doi = {10.48550/ARXIV.MATH/0201045},
}

@article{kim_stable_2019,
	title = {Stable subgroups and {Morse} subgroups in mapping class groups},
	volume = {29},
	issn = {0218-1967},
	url = {https://doi.org/10.1142/S0218196719500346},
	doi = {10.1142/S0218196719500346},
	number = {5},
	journal = {International Journal of Algebra and Computation},
	author = {Kim, Heejoung},
	year = {2019},
	mrnumber = {3978121},
	pages = {893--903},
}

@article{li_no_2020,
	title = {No growth-gaps for special cube complexes},
	volume = {14},
	issn = {1661-7207},
	url = {https://doi.org/10.4171/ggd/537},
	doi = {10.4171/ggd/537},
	number = {1},
	journal = {Groups, Geometry, and Dynamics},
	author = {Li, Jiakai and Wise, Daniel T.},
	year = {2020},
	mrnumber = {4077657},
	pages = {117--135},
}

@article{masur_geometry_1999,
	title = {Geometry of the complex of curves. {I}. {Hyperbolicity}},
	volume = {138},
	issn = {0020-9910},
	url = {https://doi.org/10.1007/s002220050343},
	doi = {10.1007/s002220050343},
	number = {1},
	journal = {Inventiones Mathematicae},
	author = {Masur, Howard A. and Minsky, Yair N.},
	year = {1999},
	mrnumber = {1714338},
	pages = {103--149},
}

@article{minsky_quasi-projections_1996,
	title = {Quasi-projections in {Teichmüller} space},
	volume = {473},
	issn = {0075-4102},
	url = {https://doi.org/10.1515/crll.1995.473.121},
	doi = {10.1515/crll.1995.473.121},
	journal = {Journal für die Reine und Angewandte Mathematik. [Crelle's Journal]},
	author = {Minsky, Yair N.},
	year = {1996},
	mrnumber = {1390685},
	pages = {121--136},
}

@article{osin_elementary_2006,
	title = {Elementary subgroups of relatively hyperbolic groups and bounded generation},
	volume = {16},
	issn = {0218-1967},
	url = {https://doi.org/10.1142/S0218196706002901},
	doi = {10.1142/S0218196706002901},
	number = {1},
	journal = {International Journal of Algebra and Computation},
	author = {Osin, Denis V.},
	year = {2006},
	mrnumber = {2217644},
	pages = {99--118},
}

@article{rafi_geodesics_2021,
	title = {Geodesics in the mapping class group},
	volume = {21},
	issn = {1472-2747},
	url = {https://doi.org/10.2140/agt.2021.21.2995},
	doi = {10.2140/agt.2021.21.2995},
	number = {6},
	journal = {Algebraic \& Geometric Topology},
	author = {Rafi, Kasra and Verberne, Yvon},
	year = {2021},
	mrnumber = {4344876},
	pages = {2995--3017},
}

@article{sisto_projections_2013,
	title = {Projections and relative hyperbolicity},
	volume = {59},
	issn = {0013-8584},
	url = {https://doi.org/10.4171/LEM/59-1-6},
	doi = {10.4171/LEM/59-1-6},
	number = {1-2},
	journal = {L'Enseignement Mathématique. Revue Internationale. 2e Série},
	author = {Sisto, Alessandro},
	year = {2013},
	mrnumber = {3113603},
	pages = {165--181},
}

@article{sisto_contracting_2018,
	title = {Contracting elements and random walks},
	volume = {742},
	issn = {0075-4102},
	url = {https://doi.org/10.1515/crelle-2015-0093},
	doi = {10.1515/crelle-2015-0093},
	journal = {Journal für die Reine und Angewandte Mathematik. [Crelle's Journal]},
	author = {Sisto, Alessandro},
	year = {2018},
	mrnumber = {3849623},
	pages = {79--114},
}

@article{yang_statistically_2019,
	title = {Statistically convex-cocompact actions of groups with contracting elements},
	issn = {1073-7928},
	url = {https://doi.org/10.1093/imrn/rny001},
	doi = {10.1093/imrn/rny001},
	number = {23},
	journal = {International Mathematics Research Notices. IMRN},
	author = {Yang, Wenyuan},
	year = {2019},
	mrnumber = {4039013},
	pages = {7259--7323},
}

@article{cordes_regularity_2022,
	title = {Regularity of {Morse} geodesics and growth of stable subgroups},
	volume = {15},
	issn = {1753-8416},
	url = {https://doi.org/10.1112/topo.12245},
	doi = {10.1112/topo.12245},
	number = {3},
	journal = {Journal of Topology},
	author = {Cordes, Matthew and Russell, Jacob and Spriano, Davide and Zalloum, Abdul},
	year = {2022},
	mrnumber = {4461848},
	pages = {1217--1247},
}

@article{vonseel_ends_2018,
	title = {Ends of {Schreier} graphs of hyperbolic groups},
	volume = {18},
	issn = {1472-2747},
	url = {https://doi.org/10.2140/agt.2018.18.3089},
	doi = {10.2140/agt.2018.18.3089},
	number = {5},
	journal = {Algebraic \& Geometric Topology},
	author = {Vonseel, Audrey},
	year = {2018},
	mrnumber = {3848409},
	pages = {3089--3118},
}

@article{russell_local--global_2022,
	title = {The local-to-global property for {Morse} quasi-geodesics},
	volume = {300},
	issn = {0025-5874},
	url = {https://doi.org/10.1007/s00209-021-02811-w},
	doi = {10.1007/s00209-021-02811-w},
	number = {2},
	journal = {Mathematische Zeitschrift},
	author = {Russell, Jacob and Spriano, Davide and Tran, Hung Cong},
	year = {2022},
	mrnumber = {4363788},
	pages = {1557--1602},
}

@article{cannon_combinatorial_1984,
	title = {The combinatorial structure of cocompact discrete hyperbolic groups},
	volume = {16},
	issn = {0046-5755},
	url = {https://doi.org/10.1007/BF00146825},
	doi = {10.1007/BF00146825},
	number = {2},
	journal = {Geometriae Dedicata},
	author = {Cannon, James W.},
	year = {1984},
	mrnumber = {758901},
	pages = {123--148},
}

@incollection{cannon_theory_1991,
	series = {Oxford {Sci}. {Publ}.},
	title = {The theory of negatively curved spaces and groups},
	booktitle = {Ergodic theory, symbolic dynamics, and hyperbolic spaces ({Trieste}, 1989)},
	publisher = {Oxford Univ. Press, New York},
	author = {Cannon, James W.},
	year = {1991},
	mrnumber = {1130181},
	pages = {315--369},
}

@misc{sisto_morse_2023,
	title = {Morse subsets of injective spaces are strongly contracting},
	copyright = {arXiv.org perpetual, non-exclusive license},
	url = {https://arxiv.org/abs/2208.13859},
	doi = {10.48550/ARXIV.2208.13859},
	publisher = {arXiv},
	author = {Sisto, Alessandro and Zalloum, Abdul},
	year = {2023},
}

@article{arzhantseva_negative_2019,
	title = {Negative curvature in graphical small cancellation groups},
	volume = {13},
	issn = {1661-7207},
	url = {https://doi.org/10.4171/GGD/498},
	doi = {10.4171/GGD/498},
	number = {2},
	journal = {Groups, Geometry, and Dynamics},
	author = {Arzhantseva, Goulnara N. and Cashen, Christopher H. and Gruber, Dominik and Hume, David},
	year = {2019},
	mrnumber = {3950644},
	pages = {579--632},
}

@misc{calvez_morse_2021,
	title = {Morse elements in {Garside} groups are strongly contracting},
	copyright = {arXiv.org perpetual, non-exclusive license},
	url = {https://arxiv.org/abs/2106.14826},
	publisher = {arXiv},
	author = {Calvez, Matthieu and Wiest, Bert},
	year = {2021},
	doi = {10.48550/ARXIV.2106.14826},
}

@misc{russell_convexity_2021,
	title = {Convexity in hierarchically hyperbolic spaces},
	url = {https://arxiv.org/abs/1809.09303},
	doi = {10.48550/ARXIV.1809.09303},
	author = {Russell, Jacob and Spriano, Davide and Tran, Hung Cong},
	year = {2021},
}

	%Xabier Legaspi
	\bigskip
	\noindent
	\emph{Xabier Legaspi} \\
	ICMAT, CSIC. 28049 Madrid, Spain.\\
	IRMAR, Université de Rennes 1. 35000 Rennes, France.\\
	\texttt{xabier.legaspi@icmat.es} \\
	\texttt{\url{https://xabierlegaspi.pages.math.cnrs.fr/}}

\end{document}